\documentclass[a4paper]{amsart}
\usepackage{color}
\usepackage[latin1]{inputenc}
\usepackage[T1]{fontenc}
\usepackage{amsfonts}
\usepackage{amssymb}
\usepackage{amsmath}
\usepackage{amsthm}
\usepackage{stmaryrd}
\usepackage{enumerate}
\usepackage{multirow}
\usepackage{graphicx}

\usepackage{graphicx,type1cm,eso-pic,color}
\usepackage{pstricks}
\usepackage{float}
\theoremstyle{theorem}
\newtheorem{theorem}{Theorem}
\newtheorem{lemma}{Lemma}
\newtheorem{proposition}{Proposition}
\newtheorem{corollary}{Corollary}

\theoremstyle{definition}
\newtheorem{remark}{Remark}

\begin{document}
\setlength\arraycolsep{2pt}
\title[N-W Estimator of the Transition Density for Diffusions]{Nadaraya-Watson Type Estimator of the Transition Density Function for Diffusion Processes}
\author{Nicolas MARIE$^{\dag}$}
\email{nmarie@parisnanterre.fr}
\author{Ousmane SACKO$^{\dag}$}
\email{osacko@parisnanterre.fr}
\date{}
\maketitle
\noindent
$^{\dag}$Universit\'e Paris Nanterre, Modal'X, 200 Avenue de la R\'epublique, 92001 Nanterre, France.
%


%
\begin{abstract}
This paper deals with a nonparametric Nadaraya-Watson (NW) estimator of the transition density function computed from independent continuous observations of a diffusion process. A risk bound is established on this estimator. The paper also deals with an extension of the penalized comparison to overfitting bandwidths selection method for our NW estimator. Finally, numerical experiments are provided.
\end{abstract}
\noindent
{\bf Keywords:} Stochastic differential equations; Transition density; Nadaraya-Watson estimator; PCO method.
%


%
\section{Introduction}\label{section_introduction}
Consider the stochastic differential equation
\begin{equation}\label{main_equation}
dX_t = b(X_t)dt +\sigma(X_t)dW_t
\textrm{ $;$ }t\in [0,2T],
\end{equation}
where $X_0 = x_0\in\mathbb R$, $W = (W_t)_{t\in [0,2T]}$ is a Brownian motion with $T > 0$, and $b,\sigma\in C^1(\mathbb R)$ with $b'$ and $\sigma'$ bounded. Under these conditions on $b$ and $\sigma$, Equation (\ref{main_equation}) has a unique (strong) solution $X = (X_t)_{t\in [0,2T]}$. Under additional conditions (see Section \ref{section_transition_density}), the transition density $p_t(x,.)$ is well-defined and can be interpreted as the conditional density of $X_{s + t}$ given $X_s = x$. For any $t\in (0,T]$, our paper deals with an adaptive Nadaraya-Watson (type) estimator of $p_t : (x,y)\mapsto p_t(x,y)$ computed from $N\in\mathbb N^*$ independent copies of $X$ observed on the time interval $[0,2T]$.\\

The copies-based statistical inference for stochastic differential equations, which is related to functional data analysis (see Ramsay and Silverman (2007) and Wang et al. (2016)), is an alternative to classic long-time behavior based methods (see Kutoyants (2004)), allowing to consider non-stationary models.

The projection least squares and the Nadaraya-Watson methods have been recently extended to the copies-based observation scheme for the estimation of $b$ (see Comte and Genon-Catalot (2020,2024), Denis et al. (2021), Comte and Marie (2023), Marie and Rosier (2023), etc.). In fact, theoretical and numerical results on such nonparametric estimators of the drift function have been also established for stochastic differential equations driven by a L\'evy process (see Halconruy and Marie (2023)) or the fractional Brownian motion (see Comte and Marie (2021) and Marie (2023)), and even for interacting particle systems (see Della Maestra and Hoffmann (2022)).

As for a regression function estimation, there are at least two kinds of nonparametric estimators of the transition density of a Markov process: those defined as a contrast minimizer, and those defined by the ratio of two density estimators extending the usual Nadaraya-Watson procedure. On the estimation of the transition density function from a discrete sample of a stationary Markov chain, the reader can refer to Lacour (2007,2008), dealing with minimum contrast estimators, or Sart (2014), dealing with both the aforementioned nonparametric estimation strategies. On the estimation of the transition density $p_t$ of Equation (\ref{main_equation}), Comte and Marie (2025) deals with a copies-based projection least squares estimator - that is a minimum contrast estimator - and the present paper deals with a copies-based Nadaraya-Watson estimator.\\

Consider
\begin{displaymath}
X^i :=\mathcal I(x_0,W^i)
\textrm{ $;$ }
i\in\{1,\dots,N\},
\end{displaymath}
where $\mathcal I(\cdot)$ is the It\^o map for Equation (\ref{main_equation}), and $W^1,\dots,W^N$ are $N$ independent copies of $W$. Consider also
\begin{eqnarray*}
 K_h(\cdot) & := &
 \frac{1}{h}K\left(\frac{.}{h}\right)
 \textrm{ $;$ }h\in (0,1],\\
 & &
 {\rm and}\quad
 Q_{\bf h} := K_{h_1}\otimes K_{h_2}
 \textrm{ $;$ }\mathbf h = (h_1,h_2)\in (0,1]^2,
\end{eqnarray*}
where $K :\mathbb R\rightarrow\mathbb R$ is a kernel function. For any $t\in (0,T]$, $\mathbf h = (h_1,h_2)\in (0,1]^2$ and $\ell\in (0,1]$, the Nadaraya-Watson estimator $\widehat p_{\mathbf h,\ell,t}$ of $p_t$ investigated in our paper is defined by
\begin{equation}\label{NW_estimator}
\widehat p_{\mathbf h,\ell,t}(x,y) :=
\frac{\widehat s_{\mathbf h,t}(x,y)}{\widehat f_{\ell}(x)}\mathbf 1_{\widehat f_{\ell}(x) >\frac{m}{2}}
\textrm{ $;$ }(x,y)\in\mathbb R^2,
\end{equation}
where $m\in (0,1]$,
\begin{eqnarray*}
 \widehat s_{\mathbf h,t}(x,y) & := &
 \frac{1}{N(T - t_0)}\sum_{i = 1}^{N}\int_{t_0}^{T}
 Q_{\bf h}(X_{s}^{i} - x,X_{s + t}^{i} - y)ds\\
 & &
 \hspace{3cm}{\rm and}\quad
 \widehat f_{\ell}(x) :=
 \frac{1}{N(T - t_0)}\sum_{i = 1}^{N}\int_{t_0}^{T}
 K_{\ell}(X_{s}^{i} - x)ds.
\end{eqnarray*}
These two last random functionals are estimators of $(x,y)\mapsto f(x)p_t(x,y)$ and
\begin{displaymath}
f(\cdot) :=\frac{1}{T - t_0}\int_{t_0}^{T}p_s(x_0,\cdot)ds
\quad {\rm respectively}.
\end{displaymath}
Note that observations of $X$ on $[0,2T]$ are required to compute $\widehat s_{\mathbf h,t}$ ($t\in [0,T]$) because for $s\in [t_0,T]$, $s + t\in [0,2T]$. Note also that the integrals involved in both the definitions of $\widehat s_{\mathbf h,t}$ and $\widehat p_{\ell}$ are considered on $[t_0,T]$ instead of $[0,T]$, because the Kusuoka-Stroock bounds on the derivatives of $p_t$ - required to control the bias terms of these estimators - explode when $t\rightarrow 0$ and are not integrable on $(0,T]$.\\

In this paper, risk bounds are established on $\widehat p_{\mathbf h,\ell,t}$ and on the adaptive estimator
\begin{displaymath}
\widehat p_{\widehat{\bf h},\widehat\ell,t}(x,y) =
\frac{\widehat s_{\widehat{\bf h},t}(x,y)}{\widehat f_{\widehat\ell}(x)}
\mathbf 1_{\widehat f_{\widehat\ell}(x) >\frac{m}{2}}
\textrm{ $;$ }(x,y)\in\mathbb R^2,
\end{displaymath}
where $\widehat{\bf h}$ (resp. $\widehat\ell$) is selected via a penalized comparison to overfitting (PCO) type criterion for its numerator (resp. denominator). The PCO bandwidth selection method, initially introduced by Lacour, Massart and Rivoirard in (2017) for the usual Parzen density estimator, needs to be extended to our framework because, contrary to the projection least squares estimator of $p_t$ investigated in Comte and Marie (2025), $\widehat p_t$ is not a minimum contrast estimator.\\

Assume that $\sigma(\cdot)^2 > 0$, and consider a known twice continuously differentiable function $v :\mathbb R\rightarrow\mathbb R$. In the spirit of Milstein et al. (2004), and as already suggested in Comte and Marie (2025) (see Section 2) for the projection least squares estimator of $p_t$, a possible application of our Nadaraya-Watson type method is to use
\begin{displaymath}
\widehat F_{\mathbf h,\ell}(t,x) :=
\int_{-\infty}^{\infty}v(y)\widehat p_{\mathbf h,\ell,T - t}(x,y)dy
\textrm{ $;$ }(t,x)\in [0,T - t_0)\times\mathbb R
\end{displaymath}
in order to solve - numerically - the parabolic partial differential equation
\begin{displaymath}
\frac{\partial u}{\partial t}(t,x) +
\frac{1}{2}\sigma(x)^2\frac{\partial^2u}{\partial x^2}(t,x) + b(x)\frac{\partial u}{\partial x}(t,x) = 0
\quad {\rm with}\quad u(T,x) = v(x),
\end{displaymath}
defining the generator of the solution with initial condition $x$ of Equation (\ref{main_equation}).\\

Section \ref{section_transition_density} deals with the existence and suitable controls of $p_t$ and $f$. Then, Sections \ref{section_non_adaptive} and \ref{section_PCO} respectively provide risk bounds on the Nadaraya-Watson estimator of $p_t$ and on its PCO-adaptive version. Finally, Section \ref{section_simulations} deals with a simulation study to show that our estimation method of $p_t$ works well.
\\
\\
{\bf Notations:}
\begin{itemize}
 \item The usual inner product (resp. norm) on $\mathbb L^2(\mathbb R^2)$ is denoted by $\langle .,.\rangle$ (resp. $\|.\|$). For the sake of readability, the usual inner product and the associated norm on $\mathbb L^2(\mathbb R)$ are denoted the same way.
 \item For a given density function $\delta :\mathbb R\rightarrow\mathbb R_+$, the usual norm on $\mathbb L^2(\mathbb R,\delta(x)dx)$ (resp. $\mathbb L^2(\mathbb R^2,\delta(y)dxdy)$) is denoted by $\|.\|_{\delta}$ (resp. $\|.\|_{1\otimes\delta}$).
\end{itemize}
%


%
\section{Preliminaries on the transition density}\label{section_transition_density}
In the sequel, $\sigma$ satisfies the following non-degeneracy condition:
\begin{equation}\label{non_degeneracy_condition}
\exists\alpha,A > 0 :\forall x\in\mathbb R\textrm{, }
\alpha\leqslant |\sigma(x)|\leqslant\alpha + A.
\end{equation}
The following lemma provides the required preliminary results on $p_t$ for our statistical purposes.
%


%
\begin{lemma}\label{properties_transition_density}
Under the condition (\ref{non_degeneracy_condition}) on $\sigma$, the transition density $p_t$, and the density function $f$ defined by
\begin{displaymath}
f(\cdot) =\frac{1}{T - t_0}\int_{t_0}^{T}p_s(x_0,\cdot)ds
\quad\textrm{with}\quad t_0\in [0,T),
\end{displaymath}
are well-defined. Moreover,
\begin{enumerate}
 \item There exists a positive constant $\overline{\mathfrak c}_T$, not depending on $t_0$, such that
 \begin{equation}\label{properties_transition_density_1}
 \sup_{(t,x,y)\in\mathbb E_{t_0}}p_t(x,y)\leqslant
 \frac{\overline{\mathfrak c}_T}{\sqrt{t_0}} =:\mathfrak m_p(t_0,T)
 \quad\textrm{with}\quad
 \mathbb E_{t_0} = [t_0,T]\times\mathbb R^2.
 \end{equation}
 \item There exists a positive constant $\mathfrak c_T$, not depending on $t_0$, such that
 \begin{equation}\label{properties_transition_density_2}
 \|f\|_{\infty}\leqslant
 \frac{2\mathfrak c_T}{\sqrt{T - t_0}} =:\mathfrak m_f(t_0,T).
 \end{equation}
 \item For every compact interval $I\subset\mathbb R$, there exists a positive constant $m$ such that $f(\cdot)\geqslant m$ on $I$.
\end{enumerate}
\end{lemma}
The proof of Lemma \ref{properties_transition_density}, relying on Menozzi et al. (2021), Theorem 1.2, is similar to that of Comte and Marie (2025), Proposition 3.1.
%


%
\begin{remark}\label{square_integrability_s_t}
Note that by Inequalities (\ref{properties_transition_density_1}) and (\ref{properties_transition_density_2}), $s_t : (x,y)\mapsto f(x)p_t(x,y)$ belongs to $\mathbb L^2(\mathbb R^2)$:
\begin{eqnarray}
 \int_{\mathbb R^2}s_t(x,y)^2dxdy
 & = &
 \int_{-\infty}^{\infty}f(x)^2\int_{-\infty}^{\infty}p_t(x,y)^2dydx
 \nonumber\\
 \label{square_integrability_s_t}
 & \leqslant &
 \|f\|_{\infty}\|p_t\|_{\infty}
 \int_{-\infty}^{\infty}f(x)
 \underbrace{\int_{-\infty}^{\infty}p_t(x,y)dy}_{= 1}dx
 \leqslant
 \mathfrak m_f(t_0,T)\mathfrak m_p(t_0,T) <\infty.
\end{eqnarray}
\end{remark}
%


%
\section{Non-adaptive risk bounds}\label{section_non_adaptive}
This section deals with non-adaptive risk bounds on $\widehat s_{\mathbf h,t}$, and then on the Nadaraya-Watson estimator $\widehat p_{\mathbf h,\ell,t}$. First, let us roughly show why $\widehat p_{\mathbf h,\ell,t}$ seems to be an appropriate estimator of $p_t$. On the one hand,
\begin{eqnarray*}
 \mathbb E(\widehat f_{\ell}(x)) & = &
 \frac{1}{T - t_0}\int_{t_0}^{T}\int_{-\infty}^{\infty}K_{\ell}(\xi - x)p_s(x_0,\xi)d\xi ds\\
 & = &
 (K_{\ell}\star f)(x)\xrightarrow[\ell\rightarrow 0]{}f(x).
\end{eqnarray*}
On the other hand, for every $s\in (0,T]$, the joint density of $(X_s,X_{s + t})$ is denoted by $p_{s,s + t}$, and since $X$ is a homogeneous Markov process,
\begin{equation}\label{consequence_homogeneity}
p_t(\xi,\zeta) =
\frac{p_{s,s + t}(\xi,\zeta)}{p_s(x_0,\xi)}
\textrm{ $;$ }\forall (\xi,\zeta)\in\mathbb R^2.
\end{equation}
Then,
\begin{eqnarray*}
 \mathbb E(\widehat s_{\mathbf h,t}(x,y)) & = &
 \frac{1}{T - t_0}\int_{t_0}^{T}\int_{\mathbb R^2}
 K_{h_1}(\xi - x)K_{h_2}(\zeta - y)p_{s,s + t}(\xi,\zeta)d\xi d\zeta ds\\
 & = &
 \int_{-\infty}^{\infty}K_{h_1}(\xi - x)
 \int_{-\infty}^{\infty}K_{h_2}(\zeta - y)p_t(\xi,\zeta)
 \underbrace{
 \frac{1}{T - t_0}\int_{t_0}^{T}p_s(x_0,\xi)ds}_{= f(\xi)}d\zeta d\xi\\
 & = &
 \int_{-\infty}^{\infty}K_{h_1}(\xi - x)(K_{h_2}\star p_t(\xi,\cdot))(y)f(\xi)d\xi\\
 & &
 \hspace{2.5cm}
 \xrightarrow[h_2\rightarrow 0]{}
 \int_{-\infty}^{\infty}K_{h_1}(\xi - x)p_t(\xi,y)f(\xi)d\xi
 \xrightarrow[h_1\rightarrow 0]{}
 f(x)p_t(x,y).
\end{eqnarray*}
For these reasons,
\begin{displaymath}
\widehat p_{\mathbf h,\ell,t} =\frac{\widehat s_{\mathbf h,t}}{\widehat f_{\ell}}
\textrm{ should be a suitable estimator of }p_t =\frac{s_t}{f}
\textrm{ with }s_t : (x,y)\mapsto f(x)p_t(x,y).
\end{displaymath}
%


%
\begin{remark}\label{relationship_joint_transition_densities}
Note that, by Equality (\ref{consequence_homogeneity}) and the definition of $f$, for every $(x,y)\in\mathbb R^2$,
\begin{eqnarray*}
 s_t(x,y) & = & f(x)p_t(x,y)\\
 & = &
 \frac{1}{T - t_0}\int_{t_0}^{T}p_s(x_0,x)p_t(x,y)ds =
 \frac{1}{T - t_0}\int_{t_0}^{T}p_{s,s + t}(x,y)ds.
\end{eqnarray*}
\end{remark}
Now, the following proposition provides risk bounds on $\widehat s_{\mathbf h,t}$.
%


%
\begin{proposition}\label{risk_bound_numerator_NW}
Assume that $K$ is a square-integrable, symmetric, kernel function. Then, for every $t\in (0,T]$,
\begin{equation}\label{risk_bound_numerator_NW_1}
\mathbb E(\|\widehat s_{\mathbf h,t} - s_t\|^2)\leqslant
\|s_{\mathbf h,t} - s_t\|^2 +\frac{\|K\|^4}{Nh_1h_2}
\end{equation}
and
\begin{equation}\label{risk_bound_numerator_NW_2}
\mathbb E(\|\widehat s_{\mathbf h,t} - s_t\|_{1\otimes f}^{2})
\leqslant\|s_{\mathbf h,t} - s_t\|_{1\otimes f}^{2} +
\frac{\|f\|_{\infty}\|K\|^4}{Nh_1h_2},
\end{equation}
where $s_{\mathbf h,t} := Q_{\bf h}\star s_t$.
\end{proposition}
Both bounds in Proposition \ref{risk_bound_numerator_NW} are derived from the bias-variance decomposition of the risk of $\widehat s_{\mathbf h,t}$. The variance term in (\ref{risk_bound_numerator_NW_1}) is of the same order as that in the risk bound on the usual Parzen estimator of the density of a $\mathbb R^2$-valued random variable (see Comte (2017), p. 65). So, (\ref{risk_bound_numerator_NW_1}) is not over-conservative. The variance term in (\ref{risk_bound_numerator_NW_2}) is a bit degraded with respect to that in (\ref{risk_bound_numerator_NW_1}) due to the change of $\mathbb L^2$-norm in the definition of the risk. Note that thanks to Inequality (\ref{properties_transition_density_2}), $\|f\|_{\infty}$ can be controlled by $\mathfrak m_f(t_0,T)$, which dependence in $t_0$ is explicit.\\

Under additional conditions on $t_0$, $K$, $b$ and $\sigma$, the following proposition provides a risk bound on $\widehat s_{\mathbf h,t}$ with an explicit bias term.
%


%
\begin{proposition}\label{control_bias_term_numerator_NW}
Assume that $t_0 > 0$, $T - t_0\geqslant 1$, and that $K$ is a square-integrable, symmetric, kernel function satisfying
\begin{equation}\label{control_bias_term_numerator_NW_1}
\int_{-\infty}^{\infty}|x^2K(x)|dx <\infty.
\end{equation}
If $(b,\sigma)$ is a smooth function, and if $(b,\sigma)$ and all its derivatives are bounded, then there exist two positive constants $\mathfrak c_{\ref{control_bias_term_numerator_NW}}$ and $q$, not depending on $t_0$, $\mathbf h$ and $N$, such that
\begin{equation}\label{control_bias_term_numerator_NW_2}
\sup_{t\in [t_0,T]}
\mathbb E(\|\widehat s_{\mathbf h,t} - s_t\|_{1\otimes f}^{2})
\leqslant
\mathfrak c_{\ref{control_bias_term_numerator_NW}}
\frac{h_{1}^{2} + h_{2}^{2}}{(1\wedge t_0)^{q + 1}} +
\frac{\mathfrak m_f(t_0,T)\|K\|^4}{Nh_1h_2},
\end{equation}
where $\mathfrak m_f(t_0,T)$ is defined in Inequality (\ref{properties_transition_density_2}).
\end{proposition}
%


%
\begin{remark}\label{remarks_control_bias_term_numerator_NW}
\white .\black
\begin{enumerate}
 \item Since the control of the bias term in Proposition \ref{control_bias_term_numerator_NW} relies on the Kusuoka-Stroock bounds on $(t,x,y)\mapsto p_t(x,y)$ and its derivatives (see Kusuoka and Stroock (1985), Corollary 3.25), the constant $\mathfrak c_{\ref{control_bias_term_numerator_NW}}$ is quite difficult to evaluate.
 \item By Proposition \ref{control_bias_term_numerator_NW}, the bias-variance tradeoff is reached by (the risk bound on) $\widehat s_{\mathbf h,t}$ when $h_1$ and $h_2$ are of order $N^{-1/4}$, leading to a rate of order $N^{-1/2}$.
 \item Consider $\beta\in\mathbb N^*$, and assume that $K$ is a square-integrable, symmetric, kernel function such that
 \begin{equation}\label{kernel_integrability_1}
 \int_{-\infty}^{\infty}|x^{\beta + 1}K(x)|dx <\infty
 \end{equation}
 and
 \begin{equation}\label{kernel_integrability_2}
 \int_{-\infty}^{\infty}x^{\upsilon}K(x)dx = 0
 \textrm{ $;$ }
 \forall\upsilon\in\{1,\dots,\beta\}.
 \end{equation}
 Such a kernel function exists by Comte (2017), Proposition 2.10. For the sake of simplicity, $\beta = 1$ in Proposition \ref{control_bias_term_numerator_NW} but, by following the same line, and thanks to the Taylor formula with integral remainder as in the proof of Marie and Rosier (2023), Proposition 1, one may establish that if $(b,\sigma)$ is a smooth function, and if $(b,\sigma)$ and all its derivatives are bounded, then there exist two positive constants $\overline{\mathfrak c}_{\ref{control_bias_term_numerator_NW}}$ and $\overline q$, not depending on $t_0$, $\mathbf h$ and $N$, such that
 \begin{equation}\label{risk_bound_numerator_NW_controlled_bias_term}
 \sup_{t\in [t_0,T]}
 \mathbb E(\|\widehat s_{\mathbf h,t} - s_t\|_{1\otimes f}^{2})
 \leqslant
 \overline{\mathfrak c}_{\ref{control_bias_term_numerator_NW}}
 \frac{h_{1}^{2\beta} + h_{2}^{2\beta}}{(1\wedge t_0)^{\overline q + 1}} +
 \frac{\mathfrak m_f(t_0,T)\|K\|^4}{Nh_1h_2}.
 \end{equation}
 Thus, the bias-variance tradeoff is reached by $\widehat s_{\mathbf h,t}$ when $h_1$ and $h_2$ are of order $N^{-1/(2\beta + 2)}$, leading to a rate of order $N^{-\beta/(\beta + 1)}$.
 \item First, to take $t_0$ close to $0$ doesn't degrade the variance term in the risk bounds (\ref{risk_bound_numerator_NW_1}), (\ref{risk_bound_numerator_NW_2}) and (\ref{control_bias_term_numerator_NW_2}). In fact, the variance term in the $\mathbb L^2$-risk bound (\ref{risk_bound_numerator_NW_1}) doesn't depend on $t_0$ at all, while that in the $f$-weighted $\mathbb L^2$-risk bound (\ref{control_bias_term_numerator_NW_2}) degrades when $t_0$ is close to $T$ because it depends on $\mathfrak m_f(t_0,T) = 2\mathfrak c_T(T - t_0)^{-1/2}$, which is the control of $\|f\|_{\infty}$ derived from the Kusuoka-Stroock bound on $p_t$ ($t\in (0,T]$) in Lemma \ref{properties_transition_density}.(2). Now, as already mentioned, the control of the bias term in Proposition \ref{control_bias_term_numerator_NW} relies on the Kusuoka-Stroock bounds on $(t,x,y)\mapsto p_t(x,y)$ and its derivatives, which are singular on $\{(0,x,x)\textrm{ $;$ }x\in\mathbb R\}$. Precisely, under the conditions of Proposition \ref{control_bias_term_numerator_NW}, the risk bound on $\widehat s_{\mathbf h,t}$ only depend on $t_0$ through a multiplicative constant of order $(1\wedge t_0)^{-q - 1}$. Therefore, by taking $t_0\in [1,T - 1]$ when $T > 2$, the risk bound (\ref{control_bias_term_numerator_NW_2}) on $\widehat s_{\mathbf h,t}$ doesn't depend on $t_0$.
\end{enumerate}
\end{remark}
Finally, the following proposition provides risk bounds on $\widehat p_{\mathbf h,\ell,t}$.
%


%
\begin{corollary}\label{risk_bound_NW}
Consider ${\tt l},{\tt r}\in\mathbb R$ satisfying ${\tt l} < {\tt r}$, and assume that $f(\cdot) > m > 0$ on $[{\tt l},{\tt r}]$. Assume also that $t_0 > 0$. Under the conditions of Proposition \ref{risk_bound_numerator_NW}, for every $t\in (0,T]$,
\begin{equation}\label{risk_bound_NW_1}
\mathbb E(\|\widehat p_{\mathbf h,\ell,t} - p_t\|_{[{\tt l},{\tt r}]\times\mathbb R}^{2})
\leqslant
\frac{16}{m^2}
\max\{1,\mathfrak m_p(t_0,T)\}
\left(\|s_{\mathbf h,t} - s_t\|^2 +\|f_{\ell} - f\|^2 +
\frac{\|K\|^4}{Nh_1h_2} +\frac{\|K\|^2}{N\ell}\right)
\end{equation}
and
\begin{eqnarray}
 \label{risk_bound_NW_2}
 \mathbb E(\|\widehat p_{\mathbf h,\ell,t} -
 p_t\|_{1\otimes f,[{\tt l},{\tt r}]\times\mathbb R}^{2})
 & \leqslant &
 \frac{16}{m^2}
 \max\{1,\mathfrak m_p(t_0,T)\}\mathfrak m_f(t_0,T)\\
 & &
 \hspace{0.75cm}\times
 \left(\frac{1}{m}\|s_{\mathbf h,t} - s_t\|_{1\otimes f}^{2} +\|f_{\ell} - f\|^2 +
 \frac{\|K\|^4}{Nh_1h_2} +\frac{\|K\|^2}{N\ell}\right),
 \nonumber
\end{eqnarray}
where $f_{\ell} := K_{\ell}\star f$ and $\mathfrak m_p(t_0,T)$ (resp. $\mathfrak m_f(t_0,T)$) is given in (\ref{properties_transition_density_1}) (resp. in (\ref{properties_transition_density_2})). 
\end{corollary}
Corollary \ref{risk_bound_NW} says that the risk of $\widehat p_{\mathbf h,\ell,t}$ is controlled by the sum of those of $\widehat s_{\mathbf h,t}$ and $\widehat f_{\ell}$ up to a multiplicative constant.
%


%
\begin{remark}\label{remarks_risk_bound_NW}
With the notations of Corollary \ref{risk_bound_NW},
\begin{enumerate}
 \item Consider $\beta\in\mathbb N^*$, assume that $K$ is a square-integrable, symmetric, kernel function satisfying the conditions (\ref{kernel_integrability_1}) and (\ref{kernel_integrability_2}), and assume also that $(b,\sigma)$ is a bounded smooth function with all its derivatives bounded. By Inequality (\ref{risk_bound_numerator_NW_controlled_bias_term}), by Marie and Rosier (2023), Proposition 1, and by Corollary \ref{risk_bound_NW},
 \begin{displaymath}
 \sup_{t\in [t_0,T]}
 \mathbb E(\|\widehat p_{\mathbf h,\ell,t} -
 p_t\|_{1\otimes f,[{\tt l},{\tt r}]\times\mathbb R}^{2})
 \lesssim
 h_{1}^{2\beta} + h_{2}^{2\beta} +\frac{1}{Nh_1h_2} +\ell^{2\beta} +\frac{1}{N\ell}.
 \end{displaymath}
 Thus, the bias-variance tradeoff is reached by $\widehat p_{\mathbf h,\ell,t}$ when
 \begin{displaymath}
 h_1,h_2 = O(N^{-\frac{1}{2\beta + 2}})
 \quad {\rm and}\quad
 \ell = O(N^{-\frac{1}{2\beta + 1}}),
 \end{displaymath}
 leading to a rate of order
 \begin{displaymath}
 N^{-\left[\left(\frac{\beta}{\beta + 1}\right)\wedge\left(\frac{2\beta}{2\beta + 1}\right)\right]} =
 N^{-\frac{\beta}{\beta + 1}}.
 \end{displaymath}
 \item Assume that $\ell = h_1h_2$. In that case, the variance term in the risk bound of Corollary \ref{risk_bound_NW} can be compared to that in the risk bound on the projection least squares estimator $\widetilde p_t$ of $p_t$ of Comte and Marie (2025) (see Theorem 4.2), which is a random element of $\mathcal S_{\varphi,m_1}\otimes\mathcal S_{\psi,m_2}$, where $\mathcal S_{\varphi,m_1} := {\rm span}\{\varphi_1,\dots,\varphi_{m_1}\}$, $\mathcal S_{\psi,m_2} := {\rm span}\{\psi_1,\dots,\psi_{m_2}\}$, and both $(\varphi_1,\dots,\varphi_{m_1})$ and $(\psi_1,\dots,\psi_{m_2})$ are orthonormal families of $\mathbb L^2(\mathbb R)$. Indeed, since the variance term of $\widetilde p_t$ is of order
 \begin{displaymath}
 \frac{m_1\mathcal L_{\psi}(m_2)}{N}
 \quad {\rm with}\quad
 \mathcal L_{\psi}(m_2) =\sup_{z\in\mathbb R}\sum_{j = 1}^{m_2}\psi_j(z)^2,
 \end{displaymath}
 our Nadaraya-Watson estimator is theoretically at least as good as $\widetilde p_t$ when $m_2\lesssim\mathcal L_{\psi}(m_2)$ (e.g. when $(\psi_1,\dots,\psi_{m_2})$ is the trigonometric basis or the Laguerre basis).
 \item Thanks to Lemma \ref{properties_transition_density}.(3), there exists $m > 0$ such that $f(\cdot) > m$ on $[{\tt l},{\tt r}]$. Precisely, by following the same line as in the proof of Comte and Marie (2025), Proposition 3.1.(2), one can show that there exist two constants $\mathfrak c_{T}^{\star},\mathfrak m_{T}^{\star} > 0$, not depending on $t_0$, ${\tt l}$ and ${\tt r}$, such that for every $x\in [{\tt l},{\tt r}]$,
 \begin{eqnarray*}
  f(x) & > &
  \frac{\mathfrak c_{T}^{\star}}{T - t_0}\int_{t_0}^{T}s^{-\frac{1}{2}}
  \exp\left(-\frac{\mathfrak m_{T}^{\star}}{s}(x - x_0)^2\right)ds\\
  & \geqslant &
  \frac{\mathfrak c_{T}^{\star}}{T - t_0}\exp\left(
  -\frac{\mathfrak m_{T}^{\star}}{t_0}\max\{({\tt r} - x_0)^2,({\tt l} - x_0)^2\}\right)
  \int_{t_0}^{T}s^{-\frac{1}{2}}ds\\
  & = &
  \frac{2\mathfrak c_{T}^{\star}}{\sqrt T +\sqrt{t_0}}\exp\left(
  -\frac{\mathfrak m_{T}^{\star}}{t_0}\max\{({\tt r} - x_0)^2,({\tt l} - x_0)^2\}\right) =: m^{\star}.
 \end{eqnarray*}
 This lower bound provides a possible value for $m$, which dependence in $t_0$, ${\tt l}$ and ${\tt r}$ is explicit. However, since the constants $\mathfrak c_{T}^{\star}$ and $\mathfrak m_{T}^{\star}$ are unknown, $m$ must be estimated in practice. For instance, Comte (2017), Section 4.2.2 suggests to take
 \begin{displaymath}
 \widehat m_{\ell} :=\min_{x\in [{\tt l},{\tt r}]}\widehat f_{\ell}(x).
 \end{displaymath}
 Consider $\beta\in\mathbb N^*$, assume that $K$ is a square-integrable, symmetric, kernel function satisfying the conditions (\ref{kernel_integrability_1}) and (\ref{kernel_integrability_2}), and assume also that $(b,\sigma)$ is a bounded smooth function with all its derivatives bounded. As in Marie and Rosier (2023) (see p. 9), one can also take
 \begin{displaymath}
 m = m_N =
 \mathfrak cN^{-\frac{\varepsilon}{3}\cdot\frac{\beta}{\beta + 1}}
 \xrightarrow[N\rightarrow\infty]{} 0,
 \end{displaymath}
 where $\mathfrak c > 0$ is a fixed constant and $\varepsilon\in (0,1)$ is chosen as close as possible to $0$. Since $f(\cdot) > m^{\star}$ on $[{\tt l},{\tt r}]$, there exists $N_0\in\mathbb N$ such that, for every $N > N_0$, $f(\cdot) > m_N\in (0,1]$ on $[{\tt l},{\tt r}]$. So, Inequalities (\ref{risk_bound_numerator_NW_controlled_bias_term}) and (\ref{risk_bound_NW_2}), and Marie and Rosier (2023), Proposition 1, lead to the following slightly degraded rate:
 \begin{displaymath}
 N^{-(1 -\varepsilon)\frac{\beta}{\beta + 1}}.
 \end{displaymath}
 \item Consider $n\in\mathbb N^*$, and let $(t_0,t_1,\dots,t_n)$ be the dissection of $[t_0,T]$ such that
 \begin{displaymath}
 t_j := t_0 + j\cdot\frac{T - t_0}{n}
 \textrm{ $;$ }\forall j\in\{1,\dots,n\}.
 \end{displaymath}
 In practice, one can consider the following discrete-time approximate of our Nadaraya-Watson estimator of $p_t$: for $k\in\{0,\dots,n - 1\}$ such that $t\in [t_k,t_{k + 1}]$,
 \begin{equation}\label{NW_estimator}
 \widehat p_{\mathbf h,\ell,k}^{n}(x,y) :=
 \frac{\widehat s_{\mathbf h,k}^{n}(x,y)}{
 \widehat f_{\ell}^{n}(x)}\mathbf 1_{\widehat f_{\ell}^{n}(x) >\frac{m}{2}}
 \textrm{ $;$ }(x,y)\in\mathbb R^2,
 \end{equation}
 where $m\in (0,1]$,
 \begin{eqnarray*}
  \widehat s_{\mathbf h,k}^{n}(x,y) & := &
  \frac{1}{nN}\sum_{i = 1}^{N}\sum_{j = 0}^{n - 1}
  Q_{\bf h}(X_{t_j}^{i} - x,X_{t_j + t_k}^{i} - y)\\
  & &
  \hspace{3cm}{\rm and}\quad
  \widehat f_{\ell}^{n}(x) :=
  \frac{1}{nN}\sum_{i = 1}^{N}\sum_{j = 0}^{n - 1}
  K_{\ell}(X_{t_j}^{i} - x).
 \end{eqnarray*}
 Although it's out of the scope of the present paper, one could establish a risk bound on $\widehat s_{\mathbf h,k}^{n}$ by following the same line as in Marie and Rosier (2023), Section 4.
\end{enumerate}
\end{remark}
%


%
\section{PCO bandwidths selection}\label{section_PCO}
First, let us present heuristically the penalized comparison to overfitting (PCO) method. Let $\widehat\varphi_u$ be an estimator of a function $\varphi$, depending on a parameter $u$, which needs to be selected from data in practice. As usual, the integrated $\mathbb L^2$-risk of $\widehat\varphi_u$ can be decomposed in the following way:
\begin{displaymath}
\mathbb E(\|\widehat\varphi_u -\varphi\|^2) =
\|\varphi_u -\varphi\|^2 +
\underbrace{\int {\rm var}(\widehat\varphi_u(x))dx}_{=:\mathfrak v_u}
\quad {\rm with}\quad
\varphi_u(\cdot) :=\mathbb E(\widehat\varphi_u(\cdot)).
\end{displaymath}
Moreover, assume that $\varphi_u(\cdot)\rightarrow\varphi(\cdot)$ and $\mathfrak v_u\rightarrow\infty$ when $u\rightarrow 0$, and that $\widehat\varphi_u$ reaches the bias-variance tradeoff for $u = u_{\star}$. In penalization-based selection methods, $u_{\star}$ is approximated by the minimizer - in a finite subset $\mathcal U$ whose element closest to $0$ is denoted by $u_0$ - of the sum of an estimator of the integrated squared-bias $\|\varphi_u -\varphi\|^2$ and of a penalty of same order as the integrated variance $\mathfrak v_u$. If $\widehat\varphi_u$ is defined as the minimizer of an objective function $\gamma(\cdot)$ (e.g. the projection estimators), then $\|\varphi_u -\varphi\|^2$ is estimated by $\gamma(\widehat\varphi_u)$. Otherwise, since $\varphi_u(\cdot)\rightarrow\varphi(\cdot)$ when $u\rightarrow 0$, $\|\varphi_u -\varphi\|\approx\|\varphi_u -\varphi_{u_0}\|$, and then $\|\varphi_u -\varphi\|^2$ can also be estimated by $\|\widehat\varphi_u -\widehat\varphi_{u_0}\|^2$. This strategy, called the penalized comparison to overfitting method, has been introduced by C. Lacour, P. Massart and V. Rivoirard in (2017) for the density of a finite-dimensional random variable estimation.\\

Now, let us provide an extension of the PCO method to the bandwidths selection of the numerator and denominator of our Nadaraya-Watson type estimator of $p_t$. Throughout this section, $t_0 > 0$. Let $\mathcal H_N$ be a finite subset of $[h_0,1]$, where $Nh_{0}^{2}\geqslant 1$. Moreover, consider $\mathbf h_0 = (h_0,h_0)$,
\begin{equation}\label{PCO_criterion_numerator}
\widehat{\bf h} =
\widehat{\bf h}(t) :=
\arg\min_{\mathbf h\in\mathcal H_{N}^{2}}\{
\|\widehat s_{\mathbf h,t} -\widehat s_{\mathbf h_0,t}\|^2 + {\rm pen}(\mathbf h)\}
\end{equation}
with
\begin{eqnarray}\label{PCO_penalty_numerator}
 {\rm pen}(\mathbf h) & = &
 \frac{2}{(T - t_0)^2N^2}\sum_{i = 1}^{N}\left\langle
 \int_{t_0}^{T}Q_{\bf h}(X_{s}^{i} -\cdot,X_{s + t}^{i} -\cdot)ds,\right.\\
 & &
 \hspace{4.5cm}\left.
 \int_{t_0}^{T}Q_{\mathbf h_0}(X_{s}^{i} -\cdot,X_{s + t}^{i} -\cdot)ds\right\rangle
 \textrm{ $;$ }\forall\mathbf h\in\mathcal H_{N}^{2},
 \nonumber
\end{eqnarray}
and
\begin{equation}\label{PCO_criterion_denominator}
\widehat\ell :=
\arg\min_{\ell\in\mathcal H_N}\{
\|\widehat f_{\ell} -\widehat f_{h_0}\|^2 + {\rm pen}^{\dag}(\ell)\}
\end{equation}
with
\begin{displaymath}
{\rm pen}^{\dag}(\ell) =
\frac{2}{(T - t_0)^2N^2}\sum_{i = 1}^{N}\left\langle
\int_{t_0}^{T}K_{\ell}(X_{s}^{i} -\cdot)ds,
\int_{t_0}^{T}K_{h_0}(X_{s}^{i} -\cdot)ds\right\rangle
\textrm{ $;$ }\forall\ell\in\mathcal H_N.
\end{displaymath}
In the PCO (bandwidth selection) criterion (\ref{PCO_criterion_numerator}), the {\it overfitting} loss $\mathbf h\mapsto\|\widehat s_{\mathbf h,t} -\widehat s_{\mathbf h_0,t}\|^2$, which models the risk to select $\mathbf h\in\mathcal H_{N}^{2}$ too close to $\mathbf h_0$ - and then to degrade excessively the variance of $\widehat s_{\mathbf h,t}$ - is penalized by ${\rm pen}(\mathbf h)$ which is of same order as the variance term in Proposition \ref{risk_bound_numerator_NW}. Regardless of the chosen penalty, the PCO criterion (\ref{PCO_criterion_numerator}) leads to the following control of the $\mathbb L^2$-loss of $\widehat s_{\widehat{\bf h},t}$ (see the proof of Theorem \ref{risk_bound_numerator_PCO_NW}):
\begin{displaymath}
\|\widehat s_{\widehat{\bf h},t} - s_t\|^2
\leqslant
\|\widehat s_{\mathbf h,t} - s_t\|^2 -\psi(\mathbf h) +\psi(\widehat{\bf h})
\textrm{ $;$ }\mathbf h\in\mathcal H_{N}^{2}
\end{displaymath}
with
\begin{eqnarray*}
 \psi(\cdot) & = &
 2\langle\widehat s_{.,t} - s_t,\widehat s_{\mathbf h_0,t} - s_t\rangle -
 \textrm{pen}(\cdot)\\
 & = &
 \psi_1(\cdot) +\psi_2(\cdot) +\psi_3(\cdot),
\end{eqnarray*}
where $\psi_1(\cdot)$ and $\psi_3(\cdot)$ are defined in Section \ref{steps_proof_risk_bound_numerator_PCO_NW} and don't depend on ${\rm pen}(\cdot)$, and
\begin{eqnarray*}
 \psi_2(\mathbf h) & := &
 \frac{2}{N^2}\sum_{i = 1}^{N}
 \left\langle\frac{1}{T - t_0}
 \int_{t_0}^{T}Q_{\bf h}(X_{s}^{i} -\cdot,X_{s + t}^{i} -\cdot)ds - s_{\mathbf h,t},\right.\\
 & &
 \hspace{3cm}\left.\frac{1}{T - t_0}
 \int_{t_0}^{T}Q_{\mathbf h_0}(X_{s}^{i} -\cdot,X_{s + t}^{i} -\cdot)ds - s_{\mathbf h_0,t}\right\rangle
 - {\rm pen}(\mathbf h)
 \textrm{ $;$ }\forall\mathbf h\in\mathcal H_{N}^{2}.
\end{eqnarray*}
So, (\ref{PCO_penalty_numerator}) provides an appropriate penalty, allowing to control $\psi_1$, $\psi_2$ and $\psi_3$ thanks to the weak Bernstein inequality (see Massart (2007), Proposition 2.9 and Inequality (2.23)) and to a concentration inequality for U-statistics (see Gin\'e and Nickl (2015), Theorem 3.4.8).
%


%
\begin{remark}\label{remark_condition_h_0}
Note that the condition $Nh_{0}^{2}\geqslant 1$ ensures that even in the worst overfitting case (i.e. $h_1 = h_2 = h_0$), the variance term in the risk bounds on $\widehat s_{\mathbf h,t}$ (see Propositions \ref{risk_bound_numerator_NW} and \ref{control_bias_term_numerator_NW}) doesn't explode when $N\rightarrow\infty$. This is crucial in order to establish the controls in Lemmas \ref{control_U_term}, \ref{control_V_term} and \ref{control_W_term} from Proposition \ref{Kernel_type_properties_Phi}.
\end{remark}
A risk bound on $\widehat f_{\widehat\ell}$ has been already established in Marie and Rosier (2023) (see Theorem 1). So, this section deals with risk bounds on $\widehat s_{\widehat{\bf h},t}$ (see Theorem \ref{risk_bound_numerator_PCO_NW}) and $\widehat p_{\widehat{\bf h},\widehat{\ell},t}$ (see Corollary \ref{risk_bound_PCO_NW}).\\

Recall that $s_t\in\mathbb  L^2(\mathbb R^2)$ by Inequalities (\ref{properties_transition_density_1}) and (\ref{properties_transition_density_2}).
%


%
\begin{theorem}\label{risk_bound_numerator_PCO_NW}
Assume that $K$ is a square-integrable, symmetric, kernel function. Then, there exist two positive (and deterministic) constants $\mathfrak c_{\ref{risk_bound_numerator_PCO_NW}}$ and $\mathfrak m_{\ref{risk_bound_numerator_PCO_NW}}$, not depending on $N$ and $t$ (but on $t_0$), such that for every $\theta\in (0,1)$ and $\lambda > 0$, with probability larger than $1 -\mathfrak m_{\ref{risk_bound_numerator_PCO_NW}}|\mathcal H_N|^2e^{-\lambda}$,
\begin{displaymath}
\|\widehat s_{\widehat{\bf h},t} - s_t\|^2
\leqslant
(1 +\theta)\min_{\mathbf h\in\mathcal H_{N}^{2}}\|\widehat s_{\mathbf h,t} - s_t\|^2 +
\frac{\mathfrak c_{\ref{risk_bound_numerator_PCO_NW}}}{\theta}\left(
\|s_{\mathbf h_0,t} - s_t\|^2 +\frac{(1 +\lambda)^3}{N}\right).
\end{displaymath}
\end{theorem}
%


%
\begin{remark}\label{explicit_risk_bound_numerator_PCO_NW}
Under the conditions of Proposition \ref{control_bias_term_numerator_NW}, by Theorem \ref{risk_bound_numerator_PCO_NW}, and since $h_0\leqslant\min(\mathcal H_N)$,
\begin{eqnarray*}
 \sup_{t\in [t_0,T]}
 \mathbb E(\|\widehat s_{\widehat{\bf h},t} - s_t\|_{1\otimes f}^{2})
 & \lesssim &
 \min_{\mathbf h\in\mathcal H_{N}^{2}}
 \left\{h_{1}^{2} + h_{2}^{2} +\frac{1}{Nh_1h_2}\right\} +
 \underbrace{\sup_{t\in [t_0,T]}\|s_{\mathbf h_0,t} - s_t\|^2}_{\leqslant 2h_{0}^{2}} +\frac{1}{N}\\
 & \lesssim &
 2\min_{\mathbf h\in\mathcal H_{N}^{2}}
 \left\{h_{1}^{2} + h_{2}^{2} +\frac{\|K\|^4}{Nh_1h_2}\right\} +\frac{1}{N}.
\end{eqnarray*}
Thus, the performance of the estimator $\widehat s_{\widehat{\bf h},t}$ is of same order as that of the best estimator in the collection $\{\widehat s_{\mathbf h,t}\textrm{ $;$ }\mathbf h\in\mathcal H_{N}^{2}\}$.
\end{remark}
%


%
\begin{corollary}\label{risk_bound_PCO_NW}
Consider ${\tt l},{\tt r}\in\mathbb R$ satisfying ${\tt l} < {\tt r}$, and assume that $f(\cdot) > m > 0$ on $[{\tt l},{\tt r}]$. Under the conditions of Theorem \ref{risk_bound_numerator_PCO_NW}, there exists a constant $\mathfrak c_{\ref{risk_bound_PCO_NW}} > 0$, not depending on $N$, $t$, ${\tt l}$ and ${\tt r}$, such that
\begin{eqnarray*}
 \mathbb E(\|\widehat p_{\widehat{\bf h},\widehat\ell,t} -
 p_t\|_{[{\tt l},{\tt r}]\times\mathbb R}^{2})
 & \leqslant &
 \frac{\mathfrak c_{\ref{risk_bound_PCO_NW}}}{m^2}
 \left(\min_{(\mathbf h,\ell)\in\mathcal H_{N}^{3}}\{
 \mathbb E(\|\widehat s_{\mathbf h,t} - s_t\|^2) +
 \mathbb E(\|\widehat f_{\ell} - f\|^2)\}\right.\\
 & &
 \hspace{3.5cm}\left. +
 \|s_{\mathbf h_0,t} - s_t\|^2 +
 \|f_{h_0} - f\|^2 +\frac{\log(N)^6}{N}\right).
\end{eqnarray*}
\end{corollary}
Corollary \ref{risk_bound_PCO_NW} says that the risk of $\widehat p_{\widehat{\bf h},\widehat\ell,t}$ is controlled by the sum of those of $\widehat s_{\widehat{\bf h},t}$ and $\widehat f_{\widehat\ell}$ up to a multiplicative constant.
%


%
\section{Numerical experiments}\label{section_simulations}
This section deals with a brief simulation study showing that our PCO-adaptive estimation procedure of $p_t$ works well. First, three usual models where $p_t$ can be explicitly computed are introduced in Section \ref{section_appropriate_models_simulations}, and then the numerical experiments on $\widehat p_{\widehat{\bf h},\widehat\ell,t}$ - defined by Equation (\ref{NW_estimator}) with $\widehat{\bf h}$ (resp. $\widehat\ell$) selected via the PCO criterion (\ref{PCO_criterion_numerator}) (resp. (\ref{PCO_criterion_denominator})) - are provided in Section \ref{section_implementation_results}.
%


%
\subsection{Appropriate models for numerical experiments}\label{section_appropriate_models_simulations}
Consider the $d$-dimensional Ornstein-Uhlenb-eck processes $\mathbf U^1,\dots,\mathbf U^N$, defined by
\begin{equation}\label{OU_processes}
d\mathbf U_{t}^{i} =
-\frac{r}{2}\mathbf U_{t}^{i}dt +\frac{\gamma}{2}d\mathbf W_{t}^{i}
\quad {\rm with}\quad
\mathbf U_{0}^{i}\sim\mathcal N_d\left(0,\frac{\gamma^2}{4r}\mathbf I_d\right),
\end{equation}
where $r,\gamma > 0$ and $\mathbf W^1,\dots,\mathbf W^N$ are independent $d$-dimensional Brownian motions. For $n\in\mathbb N^*$ and $\Delta > 0$, exact simulations of $\mathbf U^i$ along the dissection $\{\ell\Delta\textrm{ $;$ }\ell = 0,\dots,n\}$ of $[0,n\Delta]$ are computed via the following recursive formula:
\begin{equation}\label{exact_discretization_OU}
\mathbf U_{(j + 1)\Delta}^{i} =
e^{-\frac{r\Delta}{2}}\mathbf U_{j\Delta}^{i} +\varepsilon_{(j + 1)\Delta}^{i}
\quad {\rm with}\quad
\varepsilon_{\ell\Delta}^{i}
\sim_{\rm iid}\mathcal N_d\left(0,\frac{\gamma^2(1 - e^{-r\Delta})}{4r}\mathbf I_d\right).
\end{equation}
As in Comte and Marie (2025), we simulate discrete samples of the three following models thanks to (\ref{exact_discretization_OU}).
\begin{itemize}
 \item {\bf Model 1 (OU):} $X_{t}^{i} =\mathbf U_{t}^{i}$ with $d = 1$, $r = 2$ and $\gamma = 2$. Here, the transition density function is given by
 \begin{displaymath}
 p_{t}^{(1)}(x,y) =
 \sqrt{\frac{2r}{\pi\gamma^2(1 - e^{-rt})}}
 \exp\left(-\frac{2r}{\gamma^2(1 - e^{-rt})}
 (y - xe^{-\frac{rt}{2}})^2\right).
 \end{displaymath}
 \item {\bf Model 2:} $X_{t}^{i} =\tanh(\mathbf U_{t}^{i})$ with $d = 1$, $r = 4$ and $\gamma = 1$. Here, the transition density function is given by
 \begin{displaymath}
 p_{t}^{(2)}(x,y) =
 \frac{p_{t}^{(1)}({\rm atanh}(x),{\rm atanh}(y))}{1 - y^2}.
 \end{displaymath}
 \item {\bf Model 3 (CIR):} $X_{t}^{i} =\|\mathbf U_{t}^{i}\|_{2,d}^{2}$ with $d = 6$ and $r =\gamma = 1$. This is the so-called Cox-Ingersoll-Ross model. Here, the transition density function is given by
 \begin{displaymath}
 p_t^{(3)}(x,y)=
 c_t\exp(-c_t(xe^{-rt} + y))
 \left(\frac{y}{xe^{-rt}}\right)^{\frac{d}{4} -\frac{1}{2}}
 \mathcal I\left(\frac{d}{2} - 1,2c_t\sqrt{xye^{-rt}}\right),
 \end{displaymath}
 where
 \begin{displaymath}
 c_t :=\frac{2r}{\gamma^2(1 - e^{-rt})},
 \end{displaymath}
 and $\mathcal I(p,x)$ is the modified Bessel function of the first kind of order $p$ at point $x$ (see (20) in A\"it-Sahalia (1999)).
\end{itemize}
For all models, let us assume that $T = 10$, $t_0 = 0$, $t = 1$, $n = 500$, $\Delta = 0.02$, and that $K$ is the standard normal density function $x\mapsto (2\pi)^{-1/2}e^{-x^2/2}$, leading to
\begin{displaymath}
(K_{h_1}\star K_{h_2})(x) =
\frac{1}{\sqrt{2\pi(h_{1}^{2} + h_{2}^{2})}}\exp\left[-\frac{x^2}{2(h_{1}^{2} + h_{2}^{2})}\right]
\textrm{ $;$ }h_1,h_2 > 0.
\end{displaymath}
Then, for instance, the penalty ${\rm pen}^{\dag}(\cdot)$ involved in the PCO criterion for the estimator of $f$ is written as
\begin{eqnarray*}
 {\rm pen}^{\dag}(\ell)
 & = &
 \frac{2}{(T - t_0)^2N^2}\sum_{i = 1}^{N}\int_{t_0}^{T}\int_{t_0}^{T}
 (K_{\ell}\star K_{h_0})(X_{s}^{i} - X_{u}^{i})dsdu\\
 & = &
 \frac{2}{(T - t_0)^2N^2}\cdot\frac{1}{\sqrt{2\pi(\ell^2 + h_{0}^{2})}}
 \sum_{i = 1}^{N}\int_{t_0}^{T}\int_{t_0}^{T}
 \exp\left[-\frac{(X_{s}^{i} - X_{u}^{i})^2}{2(\ell^2 + h_{0}^{2})}\right]dsdu.
\end{eqnarray*}
Moreover, $\mathcal H_N =\{0.02k\textrm{ $;$ }k = 1,\dots,30\}$ in the sequel, and for the sake of simplicity, $\widehat{\bf h}$ is selected in $\{(h,h)\textrm{ $;$ }h\in\mathcal H_N\}\subset\mathcal H_{N}^{2}$ (isotropic case) in the PCO criterion for the estimator of $s_t : (x,y)\mapsto f(x)p_t(x,y)$.
%


%
\subsection{A brief presentation of the projection least squares estimator}\label{section_pLS_estimator}
In the next Section \ref{section_implementation_results} - for the three models presented in Section \ref{section_appropriate_models_simulations} - the MISE (Mean Integrated Squared Error) of our Nadaraya-Watson estimator is compared to that of the projection least squares estimator of $p_t$ investigated in Comte and Marie (2025), which is defined by:
\begin{displaymath}
\widehat p_{\mathbf m,t} =
\sum_{j = 1}^{m_1}\sum_{\ell = 1}^{m_2}
[\widehat\Theta_{\mathbf m,t}]_{j,\ell}(\varphi_j\otimes\psi_{\ell})
\quad {\rm with}\quad
\widehat\Theta_{\mathbf m,t} =
\widehat\Psi_{m_1}^{-1}\widehat Z_{\mathbf m,t},
\end{displaymath}
where $\mathbf m = (m_1,m_2)$ belongs to $\{1,\dots,N\}^2$, $(\varphi_1,\dots,\varphi_{m_1})$ (resp. $(\psi_1,\dots,\psi_{m_2})$) is a $m_1$-dimensional (resp. $m_2$-dimensional) orthonormal family in $\mathbb L^2(I)$ (resp. $\mathbb L^2(J)$), $I$ and $J$ are intervals of $\mathbb R$,
\begin{displaymath}
\widehat\Psi_{m_1} :=
\left(
\frac{1}{NT}\sum_{i = 1}^{N}
\int_{t_0 = 0}^{T}\varphi_j(X_{s}^{i})\varphi_{j'}(X_{s}^{i})ds\right)_{j,j'\in\{1,\dots,m_1\}}
\end{displaymath}
and
\begin{displaymath}
\widehat Z_{\mathbf m,t} :=
\left(
\frac{1}{NT}\sum_{i = 1}^{N}
\int_{t_0 = 0}^{T}\varphi_j(X_{s}^{i})\psi_{\ell}(X_{s + t}^{i})ds
\right)_{(j,\ell)\in\{1,\dots,m_1\}\times\{1,\dots,m_2\}}.
\end{displaymath}
In the sequel, $I = J = \mathbb R$, $(\varphi_j)_j$ (resp. $(\psi_{\ell})_{\ell}$) is the ($m_1$-dimensional) (resp. $m_2$-dimensional) Hermite basis, and $\mathbf m$ is selected from data thanks to the model selection procedure provided in Comte and Marie (2025), Section 5.
%


%
\subsection{Implementation and results}\label{section_implementation_results}
In our simulation study, all integrals with respect to time are approximated by Riemann sums along dissections of constant mesh of $[t_0,T]$ containing $n = 500$ points. Moreover, the MISE of our PCO-adaptive Nadaraya-Watson (PCO-NW) estimator of $p_t$ is approximated via Riemann sums along the dissection $\{x_j\textrm{ $;$ }j = 1,\dots,M\}$ ($M = 100$) of constant mesh of random intervals whose bounds depend on quantiles of the $X_{t}^{i}$'s and of the $X_{t + 1}^{i}$'s. Precisely, the MISE is computed by averaging, from 200 samples of $N$ copies of $X$, the approximation
\begin{displaymath}
\frac{\rm DXY}{M^2}
\sum_{j = 1}^{M}\sum_{k = 1}^{M}
(p_{\widehat{\bf h},\widehat\ell,t}(x_j,x_k) - p_t(x_j,x_k))^2
\quad {\rm of}\quad\|\widehat p_{\widehat{\bf h},\widehat\ell,t} - p_t\|^2,
\end{displaymath}
where ${\rm DXY} := (b\textrm X - a\textrm X)(b\textrm Y - a\textrm Y)$, $b\textrm X$ (resp. $a\textrm X$) is the 98\% (resp. 2\%) quantile of the $X_{t}^{i}$'s, and $b\textrm Y$ (resp. $a\textrm Y$) is the 99\% (resp. 1\%) quantile of the $X_{t + 1}^{i}$'s.\\

\begin{figure}[h!]
\begin{tabular}{cccc}
 & \includegraphics[width=0.47\textwidth]{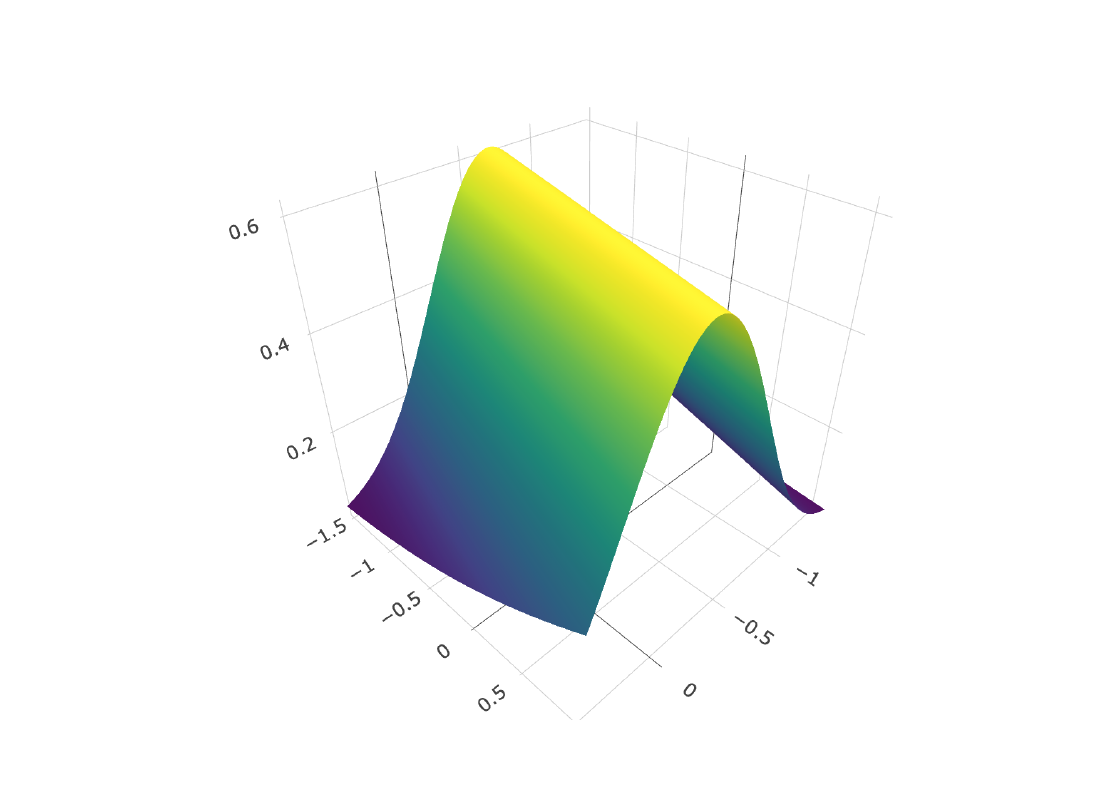} & \includegraphics[width=0.47\textwidth]{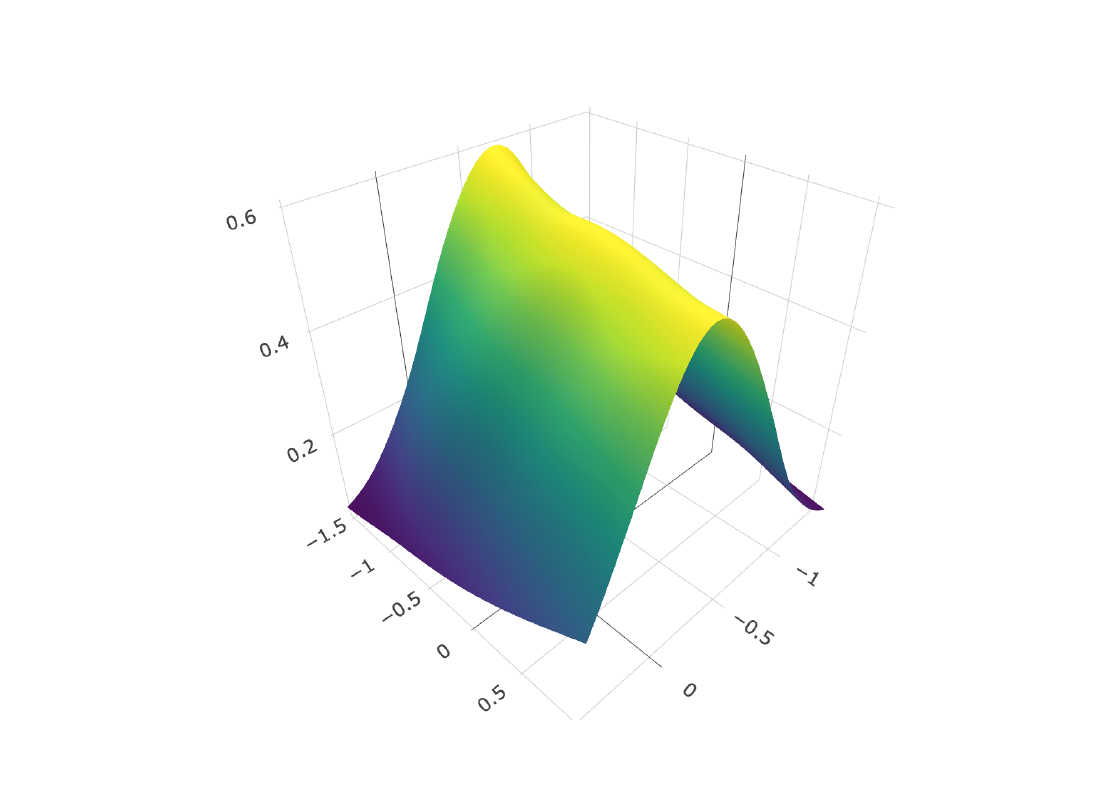}
\end{tabular}
\caption{True transition density (left) and its estimation (right) for Model 1 ($N = 200$ copies). Selected bandwidths: $\widehat h = 0.22$ and $\widehat\ell = 0.22$.}
\label{plot_model_1}
\end{figure}

\begin{figure}[h!]
\begin{tabular}{cccc}
 & \includegraphics[width=0.47\textwidth]{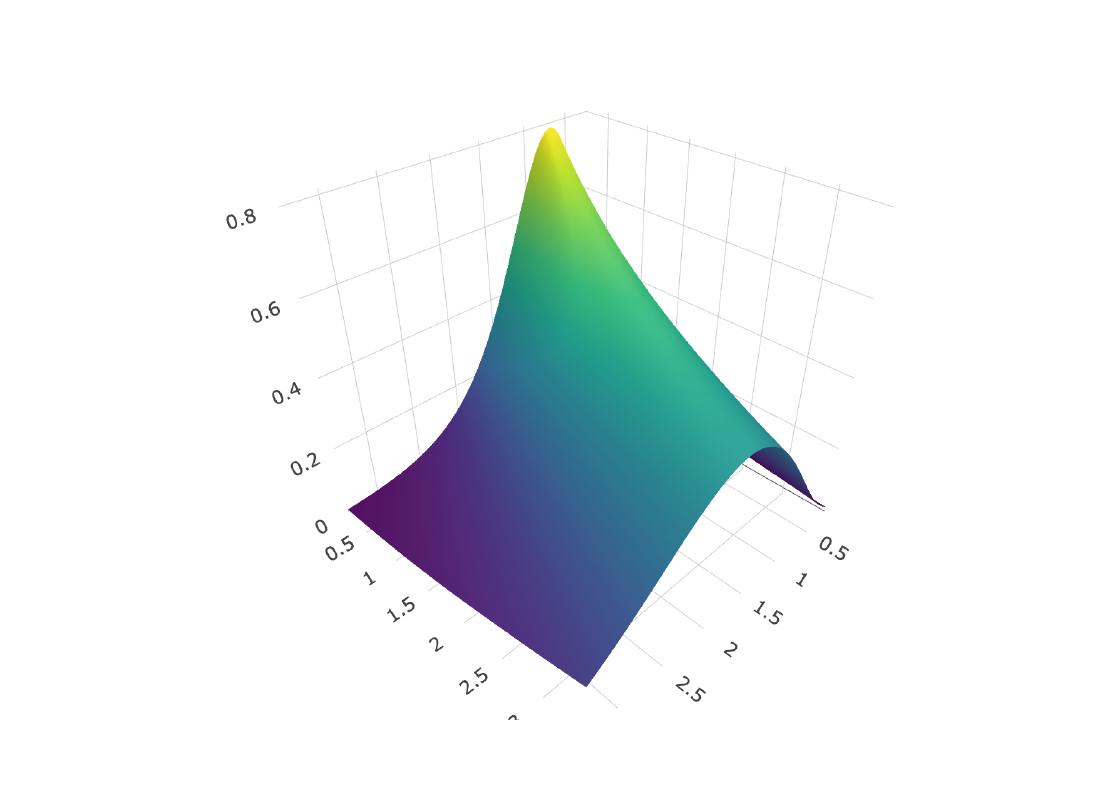} & \includegraphics[width=0.47\textwidth]{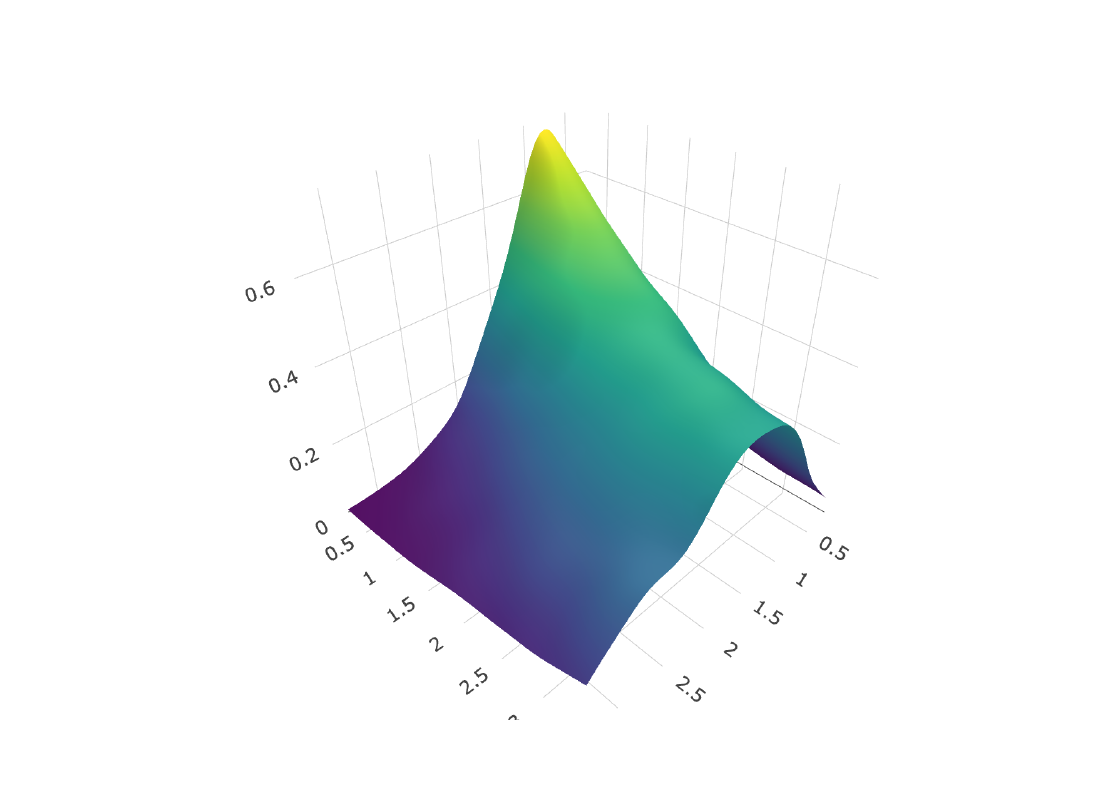}
\end{tabular}
\caption{True transition density (left) and its estimation (right) for Model 3 ($N = 200$ copies). Selected bandwidths: $\widehat h = 0.18$ and $\widehat\ell = 0.18$.}
\label{plot_model_3} 
\end{figure}

For Models 1 and 3, Figures \ref{plot_model_1} and \ref{plot_model_3} respectively display the true transition density $p_t$ on the left and its estimate obtained thanks to our procedure on the right. These figures graphically show that our PCO-NW estimator of $p_t$ works well.\\

The results of a first set of numerical experiments, involving only our PCO-NW estimator, are gathered in Table \ref{table_results_NW}. Precisely, for each model and $N\in\{100,400,1000\}$, the first line provides the MISEs - with standard deviation in parentheses - of the PCO-NW estimation of $p_t$, the second line provides the median errors, and the third (resp. fourth) line provides the mean of $\widehat h$ (resp. $\widehat\ell$).\\
For all models, the MISE is small (of order $10^{-1}$) and decreases as $N$ increases. This was expected from Corollary \ref{risk_bound_PCO_NW}. Moreover, in each situation, the median error is of same order as the MISE, illustrating the stability of our estimation procedure of $p_t$. Note also that, for all values of $N$, the MISE and the median error for Model 2 are significantly smaller than for Models 1 and 3. Finally, for all models, the means of $\widehat h$ and $\widehat\ell$ decrease as $N$ increases, but seem to stabilize from $N = 400$.\\

For each model, Table \ref{table_results_NW_vs_pLS} allows to compare - numerically - our PCO-NW estimator to the adaptive projection least squares (pLS) estimator investigated in Comte and Marie (2025) (see Section \ref{section_pLS_estimator}). When $N = 200$, the PCO-NW estimator has a smaller MISE than the pLS one for Model 1, but - in contrast - the pLS estimator performs slightly better for Model 2. When $N = 500$, except for Model 3, both estimators have similar performances, which illustrates that they are asymptotically equivalent (see Remark \ref{remarks_risk_bound_NW}.(2)). For Model 3, the PCO-NW estimator outperforms the pLS one, but probably because ${\rm supp}(p_t) =\mathbb R_{+}^{2}$ while the Hermite basis is $\mathbb R$-supported.
\begin{table}[!h]
\centering
\begin{tabular}{cl|ccc}
 Model & & $N = 100$ & $N = 400$ & $N = 1000$\\
 \hline
  & & & & \\
 1 & MISE & $0.979_{(0.535)}$ & $0.244_{(0.129)}$ & $0.140_{(0.107)}$\\
  & Medians & 0.880 & 0.223 & 0.107\\
  & Mean of $\widehat h$ & 0.275 & 0.180 & 0.120\\
  & Mean of $\widehat\ell$ & 0.288 & 0.165 & 0.140\\
  & & & & \\
 2 & MISE & $0.479_{(0.191)}$ & $0.116_{(0.057)}$ & $0.055_{(0.030)}$\\
  & Medians & 0.458 & 0.102 & 0.048\\
  & Mean of $\widehat h$ & 0.098 & 0.06 & 0.042\\
  & Mean of $\widehat\ell$ & 0.097 & 0.06 & 0.04\\
  & & & & \\
 3 & MISE & $1.052_{(2.374)}$ & $0.343_{(0.506)}$ & $0.177_{(0.414)}$\\
  & Medians & 0.771 & 0.243 & 0.102\\
  & Mean of $\widehat h$ & 0.251 & 0.127 & 0.164\\
  & Mean of $\widehat\ell$ & 0.262 & 0.107 & 0.1468\\
  & & & &
\end{tabular}
\caption{First line: 100*MISE (with 100*StD) computed over 200 repetitions. Second line: median errors. Third (resp. fourth) line: mean of $\widehat h$ (resp. $\widehat\ell$).}
\label{table_results_NW}
\end{table}
\begin{table}[!h]
\centering
\begin{tabular}{cl|cc|cc}
 Model & & \multicolumn{2}{c|}{PCO-NW} & \multicolumn{2}{c}{pLS}\\
 & & $N = 200$ & $N = 500$ & $N = 200$ & $N = 500$\\
 \hline
  & & & & \\
 1 & MISE & $0.504_{(0.349)}$ & $0.244_{(0.129)}$ & $ 0.923_{(1.606)}$ & $0.222_{(0.354)}$\\
  & Medians & 0.390 & 0.182 & 0.373 & 0.120\\
  & & & & \\
 2 & MISE & $0.282_{(0.170)}$ & $0.136_{(0.113)}$ & $0.183_{(0.211)}$ & $ 0.149_{(0.114)}$\\
  & Medians & 0.252 & 0.103 & 0.146 & 0.114\\
  & & & & \\
 3 & MISE & $0.462_{(0.276)}$ & $0.340_{(0.540)}$ & $4.067_{(8.345)}$ & $ 2.749_{(7.423)}$\\
  & Medians & 0.398 & 0.196 & 0.893 & 0.413\\
  & & & &
\end{tabular}
\caption{First line: 100*MISE (with 100*StD) computed over 200 repetitions (left: NW estimator; right: pLS estimator). Second line: median errors.}
\label{table_results_NW_vs_pLS}
\end{table}
\appendix
%


%
\section{Proofs}\label{section_proofs}
%


%
\subsection{Proof of Proposition \ref{risk_bound_numerator_NW}}
Consider $t\in (0,T]$. First of all,
\begin{displaymath}
\mathbb E(\|\widehat s_{\mathbf h,t} - s_t\|^2) =
\int_{\mathbb R^2}\mathfrak b(x,y)^2dxdy
+\int_{\mathbb R^2}\mathfrak v(x,y)dxdy
\end{displaymath}
where, for any $(x,y)\in\mathbb R^2$,
\begin{displaymath}
\mathfrak b(x,y) :=
\mathbb E(\widehat s_{\mathbf h,t}(x,y)) - s_t(x,y)
\quad {\rm and}\quad
\mathfrak v(x,y) :=
{\rm var}(\widehat s_{\mathbf h,t}(x,y)).
\end{displaymath}
First, let us find a suitable bound on the integrated squared-bias of $\widehat s_{\mathbf h,t}$. Since $X^1,\dots,X^N$ are i.i.d. copies of $X$, and by Equality (\ref{consequence_homogeneity}),
\begin{eqnarray*}
 \mathfrak b(x,y) + s_t(x,y)
 & = &
 \frac{1}{T - t_0}
 \int_{t_0}^{T}\mathbb E(Q_{\bf h}(X_s - x,X_{s + t} - y))ds\\
 & = &
 \int_{\mathbb R^2}Q_{\bf h}(\xi - x,\zeta - y)\\
 & &
 \hspace{2cm}\times
 \underbrace{\frac{1}{T - t_0}\int_{t_0}^{T}p_{s,s + t}(\xi,\zeta)ds}_{= p_t(\xi,\zeta)f(\xi)}d\xi d\zeta\\
 & = &
 \int_{\mathbb R^2}Q_{\bf h}(\xi - x,\zeta - y)s_t(\xi,\zeta)d\xi d\zeta =
 (Q_{\bf h}\star s_t)(x,y).
\end{eqnarray*}
Thus,
\begin{displaymath}
\int_{\mathbb R^2}\mathfrak b(x,y)^2dxdy =
\|s_{\mathbf h,t} - s_t\|^2
\quad {\rm with}\quad
s_{\mathbf h,t} = Q_{\bf h}\star s_t.
\end{displaymath}
Now, let us find a suitable bound on the integrated variance of $\widehat s_{\mathbf h,t}$. Again, since $X^1,\dots,X^N$ are i.i.d. copies of $X$, by Cauchy-Schwarz's (or Jensen's) inequality,
\begin{eqnarray*}
 \mathfrak v(x,y)
 & = &
 \frac{1}{N(T - t_0)^2}
 {\rm var}\left(\int_{t_0}^{T}Q_{\bf h}(X_s - x,X_{s + t} - y)ds\right)\\
 & \leqslant &
 \frac{1}{N(T - t_0)}\int_{t_0}^{T}\mathbb E(
 Q_{\bf h}(X_s - x,X_{s + t} - y)^2)ds\\
 & = &
 \frac{1}{N(T - t_0)}\int_{t_0}^{T}\left[
 \frac{1}{h_{1}^{2}h_{2}^{2}}\int_{\mathbb R^2}
 K\left(\frac{\xi - x}{h_1}\right)^2K\left(\frac{\zeta - y}{h_2}\right)^2p_{s,s + t}(\xi,\zeta)d\xi d\zeta
 \right]ds.
\end{eqnarray*}
Thus, by the Fubini-Tonelli theorem and the change of variables formula,
\begin{eqnarray*}
 \int_{\mathbb R^2}\mathfrak v(x,y)dxdy
 & = &
 \frac{1}{N(T - t_0)h_{1}^{2}h_{2}^{2}}
 \int_{\mathbb R^2}\left[\int_{\mathbb R^2}
 K\left(\frac{\xi - x}{h_1}\right)^2K\left(\frac{\zeta - y}{h_2}\right)^2dxdy\right]\\
 & &
 \hspace{6cm}\times
 \left(\int_{t_0}^{T}p_{s,s + t}(\xi,\zeta)ds\right)d\xi d\zeta\\
 & = &
 \frac{\|K\|^4}{N(T - t_0)h_1h_2}
 \int_{\mathbb R^2}
 \left(\int_{t_0}^{T}p_{s,s + t}(\xi,\zeta)ds\right)d\xi d\zeta =
 \frac{\|K\|^4}{Nh_1h_2}.
\end{eqnarray*}
Therefore,
\begin{displaymath}
\mathbb E(\|\widehat s_{\mathbf h,t} - s_t\|^2)
\leqslant
\|s_{\mathbf h,t} - s_t\|^2 +\frac{\|K\|^4}{Nh_1h_2}
\end{displaymath}
and
\begin{eqnarray*}
 \mathbb E(\|\widehat s_{\mathbf h,t} - s_t\|_{1\otimes f}^{2})
 & = &
 \int_{\mathbb R^2}\mathfrak b(x,y)^2f(y)dxdy
 +\int_{\mathbb R^2}\mathfrak v(x,y)f(y)dxdy\\
 & \leqslant &
 \|s_{\mathbf h,t} - s_t\|_{1\otimes f}^{2} +
 \frac{\|f\|_{\infty}\|K\|^4}{Nh_1h_2}.
\end{eqnarray*}
This concludes the proof.\quad $\Box$
%


%
\subsection{Proof of Proposition \ref{control_bias_term_numerator_NW}}
The proof of Proposition \ref{control_bias_term_numerator_NW} relies on the following technical lemma.
%


%
\begin{lemma}\label{Nikolskii_type_condition}
Let $\delta :\mathbb R\rightarrow\mathbb R_+$ be a density function. Under the conditions of Proposition \ref{control_bias_term_numerator_NW}, there exist two positive constants $\mathfrak c_{\ref{Nikolskii_type_condition}}$ and $q$, not depending on $t_0$, such that for every $t\in [t_0,T]$ and $\theta\in\mathbb R$,
\begin{eqnarray*}
 & &
 \int_{-\infty}^{\infty}
 \|p_t(\cdot,y +\theta) - p_t(\cdot,y)\|_{\delta}^{2}dy
 \leqslant
 \frac{\mathfrak c_{\ref{Nikolskii_type_condition}}}{t_{0}^{q}}
 (\theta^2 + |\theta|^3)\\
 & &
 \hspace{3cm}
 \textrm{and}\quad
 \int_{-\infty}^{\infty}
 \|p_t(x +\theta,\cdot) - p_t(x,\cdot)\|_{\delta}^{2}dx
 \leqslant
 \frac{\mathfrak c_{\ref{Nikolskii_type_condition}}}{t_{0}^{q}}
 (\theta^2 + |\theta|^3).
\end{eqnarray*}
\end{lemma}
The proof of Lemma \ref{Nikolskii_type_condition} is postponed to Section \ref{proof_Nikolskii_type_condition}.\\

First, by the change of variables formula and the generalized Minkowski inequality (see Comte (2017), Theorem B.1),
\begin{eqnarray*}
 \|s_{\mathbf h,t} - s_t\|_{1\otimes f}^{2} & = &
 \int_{\mathbb R^2}\left(
 \int_{\mathbb R^2}K(\xi)K(\zeta)(s_t(h_1\xi + x,h_2\zeta + y) - s_t(x,y))d\xi d\zeta\right)^2f(y)dxdy\\
 & \leqslant &
 [\int_{\mathbb R^2}|K(\xi)K(\zeta)|\\
 & &
 \hspace{1cm}\times\left(
 \int_{\mathbb R^2}(p_t(h_1\xi + x,h_2\zeta + y)f(h_1\xi + x) - p_t(x,y)f(x))^2f(y)dxdy
 \right)^{\frac{1}{2}}d\xi d\zeta]^2\\
 & \leqslant &
 3\left(\int_{\mathbb R^2}|K(\xi)K(\zeta)|
 (\mathfrak b_1(\xi,\zeta) +\mathfrak b_2(\xi,\zeta) +\mathfrak b_3(\xi,\zeta))^{\frac{1}{2}}d\xi d\zeta\right)^2
\end{eqnarray*}
where, for any $(\xi,\zeta)\in\mathbb R^2$,
\begin{eqnarray*}
 & &
 \mathfrak b_1(\xi,\zeta) :=
 \int_{\mathbb R^2}f(y)p_t(h_1\xi + x,h_2\zeta + y)^2(f(h_1\xi + x) - f(x))^2dxdy,\\
 & &
 \hspace{2cm}
 \mathfrak b_2(\xi,\zeta) :=
 \int_{\mathbb R^2}f(y)f(x)^2(p_t(h_1\xi + x,h_2\zeta + y) - p_t(x,h_2\zeta + y))^2dxdy
 \quad {\rm and}\\
 & &
 \hspace{4cm}
 \mathfrak b_3(\xi,\zeta) :=
 \int_{\mathbb R^2}f(y)f(x)^2(p_t(x,h_2\zeta + y) - p_t(x,y))^2dxdy.
\end{eqnarray*}
Now, let us find suitable bounds on $\mathfrak b_1(\xi,\zeta)$, $\mathfrak b_2(\xi,\zeta)$ and $\mathfrak b_3(\xi,\zeta)$.
\begin{itemize}
 \item\textbf{Bound on $\mathfrak b_1(\xi,\zeta)$.} By Kusuoka and Stroock (1985), Corollary 3.25, there exist two positive constants $\mathfrak c_1$ and $\mathfrak m_1$ such that, for every $t\in (0,T]$ and $(x,y)\in\mathbb R^2$,
 \begin{equation}\label{control_bias_term_numerator_NW_3}
 0 < p_t(x,y)
 \leqslant
 \frac{\mathfrak c_1}{\sqrt t}
 \exp\left[-\mathfrak m_1\frac{(y - x)^2}{t}\right].
 \end{equation}
 Then, by Inequalities (\ref{properties_transition_density_2}) and (\ref{control_bias_term_numerator_NW_3}), and since $T - t_0\geqslant 1$, for any $(\xi,\zeta)\in\mathbb R^2$,
 \begin{eqnarray*}
  \int_{-\infty}^{\infty}
  f(y)p_t(h_1\xi + x,h_2\zeta + y)^2dy
  & \leqslant &
  \|f\|_{\infty}
  \frac{\mathfrak c_{1}^{2}}{t}
  \int_{-\infty}^{\infty}\exp\left[-2\mathfrak m_1\frac{(h_2\zeta + y - h_1\xi - x)^2}{t}\right]dy\\
  & \leqslant &
  \frac{2\mathfrak c_{1}^{3}}{t}
  \int_{-\infty}^{\infty}\exp\left(-2\mathfrak m_1\frac{y^2}{t}\right)dy\\
  & \leqslant &
  \frac{\mathfrak c_2}{t_0}
  \quad {\rm with}\quad
  \mathfrak c_2 =
  2\mathfrak c_{1}^{3}
  \int_{-\infty}^{\infty}\exp\left(-2\mathfrak m_1\frac{y^2}{T}\right)dy.
 \end{eqnarray*}
 So, by Marie and Rosier (2023), Corollary 1,
 \begin{displaymath}
 \mathfrak b_1(\xi,\zeta)
 \leqslant
 \frac{\mathfrak c_2}{t_0}
 \int_{-\infty}^{\infty}
 (f(h_1\xi + x) - f(x))^2dx\\
 \leqslant
 \frac{\mathfrak c_3}{t_{0}^{r + 1}}h_{1}^{2}(\xi^2 + |\xi|^3),
 \end{displaymath}
 where $\mathfrak c_3$ and $r$ are positive constants not depending on $t_0$, $h_1$ and $\xi$.
 \item\textbf{Bounds on $\mathfrak b_2(\xi,\zeta)$ and $\mathfrak b_3(\xi,\zeta)$.} By Lemma \ref{Nikolskii_type_condition} and Inequality (\ref{properties_transition_density_2}),
 \begin{eqnarray*}
  \mathfrak b_2(\xi,\zeta) +
  \mathfrak b_3(\xi,\zeta)
  & \leqslant &
  \|f\|_{\infty}^{2}\int_{\mathbb R^2}
  f(y - h_2\zeta)(p_t(h_1\xi + x,y) - p_t(x,y))^2dxdy\\
  & &
  \hspace{2cm} +
  \|f\|_{\infty}^{2}\int_{\mathbb R^2}
  f(x)(p_t(x,h_2\zeta + y) - p_t(x,y))^2dxdy\\
  & \leqslant &
  \frac{\mathfrak c_4}{t_{0}^{q}}
  (h_{1}^{2}(\xi^2 + |\xi|^3) + h_{2}^{2}(\zeta^2 + |\zeta|^3)),
 \end{eqnarray*}
 where $\mathfrak c_4$ and $q$ are positive constants not depending on $t_0$, $\mathbf h$, $\xi$ and $\zeta$.
\end{itemize}
Thus, since $K$ is a square-integrable, symmetric, kernel function satisfying (\ref{control_bias_term_numerator_NW_1}),
\begin{eqnarray*}
 \|s_{\mathbf h,t} - s_t\|_{1\otimes f}^{2}
 & \leqslant &
 \frac{\mathfrak c_3\vee\mathfrak c_4}{t_{0}^{r + 1}\wedge t_{0}^{q}}\left(
 2h_1\int_{\mathbb R^2}|K(\xi)K(\zeta)|(\xi^2 + |\xi|^3)^{\frac{1}{2}}d\xi d\zeta\right.\\
 & &
 \hspace{4cm}\left. +
 h_2\int_{\mathbb R^2}|K(\xi)K(\zeta)|(\zeta^2 + |\zeta|^3)^{\frac{1}{2}}d\xi d\zeta\right)^2\\
 & \leqslant &
 \frac{2(\mathfrak c_3\vee\mathfrak c_4)}{(1\wedge t_0)^{q + r + 1}}(4h_{1}^{2} + h_{2}^{2})\|K\|_{1}^{2}
 \underbrace{\left(\int_{-\infty}^{\infty}|K(\xi)|(\xi^2 + |\xi|^3)^{\frac{1}{2}}d\xi\right)^2}_{<\infty}\\
 & \leqslant &
 \frac{\mathfrak c_5}{(1\wedge t_0)^{q + r + 1}}(h_{1}^{2} + h_{2}^{2}),
\end{eqnarray*}
where $\mathfrak c_5$ is a positive constant not depending on $t_0$ and $\mathbf h$. This concludes the proof.\quad $\Box$
%


%
\subsubsection{Proof of Lemma \ref{Nikolskii_type_condition}}\label{proof_Nikolskii_type_condition}
The proof of Lemma \ref{Nikolskii_type_condition} is similar to that of Marie and Rosier (2023), Corollary 1. By Kusuoka and Stroock (1985), Corollary 3.25, there exist three positive constants $\mathfrak c_1$, $\mathfrak m_1$ and $r$ such that, for every $t\in (0,T]$ and $(x,y)\in\mathbb R^2$,
\begin{equation}\label{Nikolskii_type_condition_1}
|\partial_1p_t(x,y)| +
|\partial_2p_t(x,y)|
\leqslant
\frac{\mathfrak c_1}{t^r}
\exp\left[-\mathfrak m_1\frac{(y - x)^2}{t}\right].
\end{equation}
For any $t\in [t_0,T]$ and $\vartheta\in\mathbb R_+$, by Inequality (\ref{Nikolskii_type_condition_1}),
\begin{eqnarray*}
 \int_{-\infty}^{\infty}
 \|p_t(\cdot,y +\vartheta) - p_t(\cdot,y)\|_{\delta}^{2}dy
 & \leqslant &
 \vartheta^2\int_{-\infty}^{\infty}\delta(x)\int_{-\infty}^{\infty}\left(
 \sup_{z\in [y,y +\vartheta]}|\partial_2p_t(x,x + z)|^2\right)dydx\\
 & \leqslant &
 \vartheta^2\frac{\mathfrak c_{1}^{2}}{t^{2r}}
 \left(\int_{-\infty}^{\infty}\delta(x)dx\right)
 \int_{-\infty}^{\infty}\left[
 \sup_{z\in [y,y +\vartheta]}\exp\left(-2\mathfrak m_1\frac{z^2}{t}\right)\right]dy\\
 & = &
 \vartheta^2\frac{\mathfrak c_{1}^{2}}{t^{2r}}\left[
 \int_{-\infty}^{-\vartheta}\exp\left(-2\mathfrak m_1
 \frac{(y +\vartheta)^2}{t}\right)dy +\vartheta +
 \int_{0}^{\infty}\exp\left(-2\mathfrak m_1
 \frac{y^2}{t}\right)dy\right]\\
 & \leqslant &
 \frac{\mathfrak c_{1}^{2}}{t_{0}^{2r}}(\mathfrak c_2\vartheta^2 +\vartheta^3)
\end{eqnarray*}
with
\begin{displaymath}
\mathfrak c_2 =
2\int_{0}^{\infty}\exp\left(-2\mathfrak m_1
\frac{y^2}{T}\right)dy,
\end{displaymath}
and the same way,
\begin{eqnarray*}
 \int_{-\infty}^{\infty}
 \|p_t(\cdot,y -\vartheta) - p_t(\cdot,y)\|_{\delta}^{2}dy
 & \leqslant &
 \vartheta^2\frac{\mathfrak c_{1}^{2}}{t^{2r}}\left(
 \int_{-\infty}^{0}\exp\left(-2\mathfrak m_1
 \frac{y^2}{t}\right)dy +\vartheta +
 \int_{\vartheta}^{\infty}\exp\left(-2\mathfrak m_1
 \frac{(y -\vartheta)^2}{t}\right)dy\right)\\
 & \leqslant &
 \frac{\mathfrak c_{1}^{2}}{t_{0}^{2r}}(\mathfrak c_2\vartheta^2 +\vartheta^3).
\end{eqnarray*}
Thus, for any $\theta\in\mathbb R$,
\begin{displaymath}
\int_{-\infty}^{\infty}
\|p_t(\cdot,y +\theta) - p_t(\cdot,y)\|_{\delta}^{2}dy
\leqslant
\frac{\mathfrak c_{1}^{2}}{t_{0}^{2r}}(\mathfrak c_2\theta^2 + |\theta|^3).
\end{displaymath}
By following the same line, Inequality (\ref{Nikolskii_type_condition_1}) leads to
\begin{displaymath}
\int_{-\infty}^{\infty}
\|p_t(x +\theta,\cdot) - p_t(x,\cdot)\|_{\delta}^{2}dx
\leqslant
\frac{\mathfrak c_{1}^{2}}{t_{0}^{2r}}(\mathfrak c_2\theta^2 + |\theta|^3).
\end{displaymath}
This concludes the proof.\quad $\Box$
%


%
\subsection{Proof of Corollary \ref{risk_bound_NW}}
First of all, for any $t\in (0,T]$,
\begin{displaymath}
\widehat p_{\mathbf h,\ell,t} - p_t =
\left[\frac{\widehat s_{\mathbf h,t} - s_t}{\widehat f_{\ell}} +
\left(\frac{1}{\widehat f_{\ell}} -\frac{1}{f}\right)fp_t\right]
\mathbf 1_{\widehat f_{\ell}(\cdot) >\frac{m}{2}} - 
p_t\mathbf 1_{\widehat f_{\ell}(\cdot)\leqslant\frac{m}{2}}.
\end{displaymath}
Then,
\begin{displaymath}
\|\widehat p_{\mathbf h,\ell,t} - p_t\|_{[\texttt l,\texttt r]\times\mathbb R}^{2} =
\left\|\left[\frac{\widehat s_{\mathbf h,t} - s_t}{\widehat f_{\ell}} +
\left(\frac{1}{\widehat f_{\ell}} -\frac{1}{f}\right)fp_t\right]
\mathbf 1_{\widehat f_{\ell}(\cdot) >\frac{m}{2}}\right\|_{[\texttt l,\texttt r]\times\mathbb R}^{2} +
\|p_t\mathbf 1_{\widehat f_{\ell}(\cdot)\leqslant
\frac{m}{2}}\|_{[\texttt l,\texttt r]\times\mathbb R}^{2}.
\end{displaymath}
Moreover, for any $x\in [\texttt l,\texttt r]$, since $f(x) > m$, for every $\omega\in\{\widehat f_{\ell}(\cdot)\leqslant m/2\}$,
\begin{displaymath}
|f(x) -\widehat f_{\ell}(x,\omega)|
\geqslant
f(x) -\widehat f_{\ell}(x,\omega) > m -\frac{m}{2} =\frac{m}{2}.
\end{displaymath}
Thus,
\begin{eqnarray*}
 \|\widehat p_{\mathbf h,\ell,t} - p_t\|_{[\texttt l,\texttt r]\times\mathbb R}^{2}
 & \leqslant &
 \frac{8}{m^2}\|\widehat s_{\mathbf h,t} - s_t\|_{[\texttt l,\texttt r]\times\mathbb R}^{2} +
 \frac{8}{m^2}\|(f -\widehat f_{\ell})p_t\|_{[\texttt l,\texttt r]\times\mathbb R}^{2} +
 2\|p_t\mathbf 1_{|f(\cdot) -
 \widehat f_{\ell}(\cdot)| >\frac{m}{2}}\|_{[\texttt l,\texttt r]\times\mathbb R}^{2}\\
 & \leqslant &
 \frac{8}{m^2}\int_{[\texttt l,\texttt r]\times\mathbb R}
 (\widehat s_{\mathbf h,t} - s_t)(x,y)^2dxdy\\
 & & \hspace{1cm} + 
 \frac{8}{m^2}\int_{[\texttt l,\texttt r]\times\mathbb R}(f(x) -
 \widehat f_{\ell}(x))^2p_t(x,y)^2dxdy\\
 & & \hspace{2cm} + 
 2\int_{[\texttt l,\texttt r]\times\mathbb R}p_t(x,y)^2\mathbf 1_{|f(x) -\widehat f_{\ell}(x)| >\frac{m}{2}}dxdy.
\end{eqnarray*}
By Inequality (\ref{properties_transition_density_1}), and since $p_t(x,\cdot)$ ($x\in\mathbb R$) is a density function,
\begin{displaymath}
\|\widehat p_{\mathbf h,\ell,t} - p_t\|_{[\texttt l,\texttt r]\times\mathbb R}^{2}
\leqslant
\frac{8}{m^2}
\|\widehat s_{\mathbf h,t} - s_t\|_{[\texttt l,\texttt r]\times\mathbb R}^{2} +
\frac{8}{m^2}\mathfrak m_p(t_0,T)\|\widehat f_{\ell} - f\|^2 +
2\mathfrak m_p(t_0,T)
\int_{-\infty}^{\infty}\mathbf 1_{|f(x) -\widehat f_{\ell}(x)| >\frac{m}{2}}dx.
\end{displaymath}
Therefore, by Markov's inequality,
\begin{eqnarray}
 \mathbb E(\|\widehat p_{\mathbf h,\ell,t} - p_t\|_{[\texttt l,\texttt r]\times\mathbb R}^{2})
 & \leqslant &
 \frac{8}{m^2}\mathbb E(
 \|\widehat s_{\mathbf h,t} - s_t\|_{[\texttt l,\texttt r]\times\mathbb R}^{2}) +
 \frac{8}{m^2}\mathfrak m_p(t_0,T)
 \mathbb E(\|\widehat f_{\ell} - f\|^2)
 \nonumber\\
 \label{risk_bound_NW_3}
 & &
 \hspace{4.5cm}
 +\frac{8}{m^2}\mathfrak m_p(t_0,T)
 \int_{-\infty}^{\infty}\mathbb E((f(x) -\widehat f_{\ell}(x))^2)dx
 \nonumber\\
 & \leqslant &
 \frac{8}{m^2}\max\{1,\mathfrak m_p(t_0,T)\}
 (\mathbb E(\|\widehat s_{\mathbf h,t} - s_t\|_{[\texttt l,\texttt r]\times\mathbb R}^{2}) +
 2\mathbb E(\|\widehat f_{\ell} - f\|^2)),
\end{eqnarray}
leading to Inequality (\ref{risk_bound_NW_1}) thanks to Proposition \ref{risk_bound_numerator_NW} and Marie and Rosier (2023), Proposition 1, which allow to control
\begin{displaymath}
\mathbb E(\|\widehat s_{\mathbf h,t} - s_t\|_{[\texttt l,\texttt r]\times\mathbb R}^{2})
\quad {\rm and}\quad
\mathbb E(\|\widehat f_{\ell} - f\|^2)
\quad {\rm respectively}.
\end{displaymath}
By Inequality (\ref{properties_transition_density_2}), and since $f(\cdot) > m$ on $[\texttt l,\texttt r]$,
\begin{displaymath}
m\|.\|_{[\texttt l,\texttt r]\times\mathbb R}^{2}\leqslant
\|.\|_{1\otimes f,[\texttt l,\texttt r]\times\mathbb R}^{2}\leqslant
\mathfrak m_f(t_0,T)\|.\|_{[\texttt l,\texttt r]\times\mathbb R}^{2},
\end{displaymath}
and then one may also establish Inequality (\ref{risk_bound_NW_2}) thanks to Inequality (\ref{risk_bound_NW_3}).\quad $\Box$
%


%
\subsection{Proof of Theorem \ref{risk_bound_numerator_PCO_NW}}
For any $t\in [t_0,T]$ and $\mathbf h\in (0,1]^2$, consider the map $\Phi_{\mathbf h,t}$ defined on $C^0([0,T])\times\mathbb R^2$ by
\begin{displaymath}
\Phi_{\mathbf h,t}(\varphi;x,y) :=
\frac{1}{T - t_0}\int_{t_0}^{T}Q_{\bf h}(\varphi(s) - x,\varphi(s + t) - y)ds
\end{displaymath}
for every $\varphi\in C^0([0,T])$ and $(x,y)\in\mathbb R^2$. Then,
\begin{displaymath}
\widehat s_{\mathbf h,t}(\cdot) =
\frac{1}{N}\sum_{i = 1}^{N}\Phi_{\mathbf h,t}(X^i;\cdot).
\end{displaymath}
First, the following proposition shows that
\begin{displaymath}
\mathcal K_N :=
\{(\varphi,x,y)\mapsto\Phi_{\mathbf h,t}(\varphi;x,y)
\textrm{ $;$ }\mathbf h\in\mathcal H_{N}^{2}\}
\end{displaymath}
satisfies properties close to those of a kernels set in the nonparametric regression framework (see Halconruy and Marie (2020), Assumption 2.1).
%


%
\begin{proposition}\label{Kernel_type_properties_Phi}
Under the conditions of Theorem \ref{risk_bound_numerator_PCO_NW}, there exists a constant $\mathfrak m_{\Phi} > 0$, not depending on $t$, such that:
\begin{enumerate}
 \item For every $\mathbf h = (h_1,h_2)\in (0,1]^2$ and $\varphi\in C^0([0,T])$,
 \begin{displaymath}
 \|\Phi_{\mathbf h,t}(\varphi;\cdot)\|^2
 \leqslant\frac{\mathfrak m_{\Phi}}{h_1h_2}.
 \end{displaymath}
 \item For every $\mathbf h,\mathbf l\in (0,1]^2$,
 \begin{displaymath}
 \mathbb E(\langle\Phi_{\mathbf h,t}(X^1;\cdot),\Phi_{\mathbf l,t}(X^2;\cdot)\rangle^2)
 \leqslant\mathfrak m_{\Phi}\overline s_{\mathbf l,t},
 \end{displaymath}
 where
 \begin{displaymath}
 \overline s_{\mathbf l,t} :=
 \mathbb E(\|\Phi_{\mathbf l,t}(X;\cdot)\|^2)
 \end{displaymath}
 and $(X^1,X^2)$ is a pair of independent copies of $X$.
 \item For every $\mathbf h\in (0,1]^2$ and $\varphi\in\mathbb L^2(\mathbb R^2)$,
 \begin{displaymath}
 \mathbb E(\langle\Phi_{\mathbf h,t}(X;\cdot),\varphi\rangle^2)
 \leqslant\mathfrak m_{\Phi}\|\varphi\|^2.
 \end{displaymath}
 \item For every $\mathbf h,\mathbf l\in (0,1]^2$,
 \begin{displaymath}
 |\langle\Phi_{\mathbf h,t}(X;\cdot),s_{\mathbf l,t}\rangle|
 \leqslant\mathfrak m_{\Phi},
 \end{displaymath}
 where
 \begin{displaymath}
 s_{\mathbf l,t}(\cdot) =
 \mathbb E(\widehat s_{\mathbf l,t}(\cdot)) =
 (Q_{\bf l}\star s_t)(\cdot).
 \end{displaymath}
\end{enumerate}
\end{proposition}
The proof of Proposition \ref{Kernel_type_properties_Phi} is postponed to Section \ref{proof_kernel_type_properties_Phi}. Now, the three following lemmas deal with controls of the maps $U$, $V$ and $W$ involved in the proof of Theorem \ref{risk_bound_numerator_PCO_NW} (see (the next) Section \ref{steps_proof_risk_bound_numerator_PCO_NW}).
%


%
\begin{lemma}\label{control_U_term}
For every $\mathbf h,\mathbf l\in\mathcal H_{N}^{2}$, consider
\begin{displaymath}
U_{\mathbf h,\mathbf l} :=\sum_{i\neq k}
\langle\Phi_{\mathbf h,t}(X^i;\cdot) - s_{\mathbf h,t},
\Phi_{\mathbf l,t}(X^k;\cdot) - s_{\mathbf l,t}\rangle.
\end{displaymath}
There exists a deterministic constant $\mathfrak c_{\ref{control_U_term}} > 0$, not depending on $N$ and $t$, such that for every $\theta\in (0,1)$ and $\lambda > 0$, with probability larger than $1 - 5.4|\mathcal H_N|^2e^{-\lambda}$,
\begin{eqnarray*}
 \sup_{\mathbf h\in\mathcal H_{N}^{2}}
 \left\{\frac{|U_{\mathbf h,\mathbf h_0}|}{N^2}
 -\frac{\theta}{N}\overline s_{\mathbf h,t}\right\}
 & \leqslant &
 \frac{\mathfrak c_{\ref{control_U_term}}(1 +\lambda)^3}{\theta N}\\
 & &
 \textrm{and}\quad
 \sup_{\mathbf h\in\mathcal H_{N}^{2}}
 \left\{\frac{|U_{\mathbf h,\mathbf h}|}{N^2}
 -\frac{\theta}{N}\overline s_{\mathbf h,t}\right\}
 \leqslant
 \frac{\mathfrak c_{\ref{control_U_term}}(1 +\lambda)^3}{\theta N}.
\end{eqnarray*}
\end{lemma}
%


%
\begin{lemma}\label{control_V_term}
For every $\mathbf h\in\mathcal H_{N}^{2}$, consider
\begin{displaymath}
V_{\bf h} :=\frac{1}{N}\sum_{i = 1}^{N}
\|\Phi_{\mathbf h,t}(X^i;\cdot) - s_{\mathbf h,t}\|^2.
\end{displaymath}
There exists a deterministic constant $\mathfrak c_{\ref{control_V_term}} > 0$, not depending on $N$ and $t$, such that for every $\theta\in (0,1)$ and $\lambda > 0$, with probability larger than $1 - 2|\mathcal H_N|^2e^{-\lambda}$,
\begin{displaymath}
\sup_{\mathbf h\in\mathcal H_{N}^{2}}\left\{
\frac{1}{N}|V_{\bf h} -\overline s_{\mathbf h,t}| -\frac{\theta}{N}\overline s_{\mathbf h,t}\right\}
\leqslant
\frac{\mathfrak c_{\ref{control_V_term}}(1 +\lambda)}{\theta N}.
\end{displaymath}
\end{lemma}
%


%
\begin{lemma}\label{control_W_term}
For every $\mathbf h,\mathbf l\in\mathcal H_{N}^{2}$, consider
\begin{displaymath}
W_{\mathbf h,\mathbf l} :=
\langle\widehat s_{\mathbf h,t} - s_{\mathbf h,t},s_{\mathbf l,t} - s_t\rangle.
\end{displaymath}
There exists a deterministic constant $\mathfrak c_{\ref{control_W_term}} > 0$, not depending on $N$ and $t$, such that for every $\theta\in (0,1)$ and $\lambda > 0$, with probability larger than $1 - 2|\mathcal H_N|^2e^{-\lambda}$,
\begin{eqnarray*}
 & &
 \sup_{\mathbf h\in\mathcal H_{N}^{2}}\{
 |W_{\mathbf h,\mathbf h_0}| -\theta\|s_{\mathbf h_0,t} - s_t\|^2\}
 \leqslant
 \frac{\mathfrak c_{\ref{control_W_term}}(1 +\lambda)^2}{\theta N},\\
 & &
 \hspace{2cm}
 \sup_{\mathbf h\in\mathcal H_{N}^{2}}\{
 |W_{\mathbf h_0,\mathbf h}| -\theta\|s_{\mathbf h,t} - s_t\|^2\}
 \leqslant
 \frac{\mathfrak c_{\ref{control_W_term}}
 (1 +\lambda)^2}{\theta N}\\
 & &
 \hspace{4cm}{\rm and}\quad
 \sup_{\mathbf h\in\mathcal H_{N}^{2}}\{
 |W_{\mathbf h,\mathbf h}| -\theta\|s_{\mathbf h,t} - s_t\|^2\}
 \leqslant
 \frac{\mathfrak c_{\ref{control_W_term}}
 (1 +\lambda)^2}{\theta N}.
\end{eqnarray*}
\end{lemma}
As in Marie and Rosier (2023) (see Lemmas 4, 5 and 6), the proofs of Lemmas \ref{control_U_term}, \ref{control_V_term} and \ref{control_W_term} rely on Proposition \ref{Kernel_type_properties_Phi}, on a concentration inequality for U-statistics (see Gin\'e and Nickl (2015), Theorem 3.4.8), and on the weak Bernstein inequality (see Massart (2007), Proposition 2.9 and Inequality (2.23)). So, the proofs of Lemmas \ref{control_U_term}, \ref{control_V_term} and \ref{control_W_term} are omitted.
%


%
\begin{remark}\label{remarks_proof_risk_bound_numerator_PCO_NW}
\white .\black
\begin{enumerate}
 \item Precisely, in the proof of Theorem \ref{risk_bound_numerator_PCO_NW} (see Section \ref{steps_proof_risk_bound_numerator_PCO_NW}), the maps $U$, $V$ and $W$ are involved in the control of the $\mathbb L^2$-loss of $\widehat s_{\widehat{\bf h},t}$ through
 \begin{displaymath}
 \|\widehat s_{\mathbf h,t} - s_{\mathbf h,t}\|^2 =
 \frac{U_{\mathbf h,\mathbf h}}{N^2} +\frac{V_{\bf h}}{N}
 \textrm{ $;$ }\mathbf h\in\mathcal H_{N}^{2},
 \end{displaymath}
 and
 \begin{eqnarray*}
  \langle\widehat s_{\mathbf h,t} - s_t,\widehat s_{\mathbf l,t} - s_t\rangle
  & = &
  U_{\mathbf h,\mathbf l} + W_{\mathbf h,\mathbf l} + W_{\mathbf l,\mathbf h} +
  \langle s_{\mathbf h,t} - s_t,s_{\mathbf l,t} - s_t\rangle\\
  & &
  \hspace{1cm} +
  \frac{1}{N^2}\sum_{i = 1}^{N}\langle\Phi_{\mathbf h,t}(X^i;\cdot) - s_{\mathbf h,t},
  \Phi_{\mathbf l,t}(X^i;\cdot) - s_{\mathbf l,t}\rangle
  \textrm{ $;$ }\mathbf h,\mathbf l\in\mathcal H_{N}^{2}.
 \end{eqnarray*}
 \item Consider $\mathfrak m\in\mathbb R$, $\mathfrak c\in (0,1)$ and a set of random variables $\{\xi_{\bf h}\textrm{ $;$ }\mathbf h\in\mathcal H_{N}^{2}\}$ such that
 \begin{displaymath}
 \mathbb P(\xi_{\bf h}\leqslant\mathfrak m) > 1 -\mathfrak c
 \textrm{ $;$ }
 \forall\mathbf h\in\mathcal H_{N}^{2}.
 \end{displaymath}
 Then,
 \begin{eqnarray*}
  \mathbb P\left(\sup_{\mathbf h\in\mathcal H_{N}^{2}}\xi_{\mathbf h}\leqslant\mathfrak m\right)
  & = &
  1 -\mathbb P\left(
  \bigcup_{\mathbf h\in\mathcal H_{N}^{2}}\{\xi_{\mathbf h} >\mathfrak m\}\right)\\
  & \geqslant &
  1 -\sum_{\mathbf h\in\mathcal H_{N}^{2}}(1 -\mathbb P(\xi_{\mathbf h}\leqslant\mathfrak m)) >
  1 - |\mathcal H_N|^2\mathfrak c.
 \end{eqnarray*}
 In Lemmas \ref{control_U_term}, \ref{control_V_term} and \ref{control_W_term}, this is the way the bounds on supremums are derived from bounds established for $\mathbf h\in\mathcal H_{N}^{2}$ fixed.
\end{enumerate}
\end{remark}
%


%
\subsubsection{Steps of the proof of Theorem \ref{risk_bound_numerator_PCO_NW}}\label{steps_proof_risk_bound_numerator_PCO_NW}
The proof of Theorem \ref{risk_bound_numerator_PCO_NW} is dissected in four steps. Step 1 shows that, for any $\mathbf h\in\mathcal H_{N}^{2}$,
\begin{displaymath}
\|\widehat s_{\widehat{\bf h},t} - s_t\|^2
\leqslant\|\widehat s_{\mathbf h,t} - s_t\|^2 -\psi(\mathbf h) +\psi(\widehat{\bf h}),
\end{displaymath}
where $\psi$ is a map depending on $U$ and $W$. Then, $\psi(\mathbf h)$ and $\psi(\widehat{\bf h})$ are controlled in Step 2 thanks to Lemmas \ref{control_U_term} and \ref{control_W_term}. Step 3 deals with a two-sided relationship between
\begin{displaymath}
\|\widehat s_{\mathbf h,t} - s_t\|^2
\quad {\rm and}\quad
\|s_{\mathbf h,t} - s_t\|^2
\textrm{ $;$ }\mathbf h\in\mathcal H_{N}^{2},
\end{displaymath}
thanks to Lemmas \ref{control_U_term}, \ref{control_V_term} and \ref{control_W_term}. The conclusion comes in Step 4.
\\
\\
\textbf{Step 1.} First,
\begin{displaymath}
\|\widehat s_{\widehat{\bf h},t} - s_t\|^2 =
\|\widehat s_{\widehat{\bf h},t} -\widehat s_{\mathbf h_0,t}\|^2 +
\|\widehat s_{\mathbf h_0,t} - s_t\|^2 +
2\langle\widehat s_{\widehat{\bf h},t} -\widehat s_{\mathbf h_0,t},
\widehat s_{\mathbf h_0,t} - s_t\rangle
\end{displaymath}
and, for any $\mathbf h\in\mathcal H_{N}^{2}$,
\begin{eqnarray*}
 \|\widehat s_{\widehat{\bf h},t} -\widehat s_{\mathbf h_0,t}\|^2
 & \leqslant &
 \|\widehat s_{\mathbf h,t} -\widehat s_{\mathbf h_0,t}\|^2 +
 {\rm pen}(\mathbf h) - {\rm pen}(\widehat{\bf h})
 \quad\textrm{by (\ref{PCO_criterion_numerator})}\\
 & = &
 \|\widehat s_{\mathbf h,t} - s_t\|^2 +
 2\langle\widehat s_{\mathbf h,t} - s_t,s_t -\widehat s_{\mathbf h_0,t}\rangle +
 \|s_t -\widehat s_{\mathbf h_0,t}\|^2 +
 {\rm pen}(\mathbf h) - {\rm pen}(\widehat{\bf h})\\
 & = &
 \|\widehat s_{\mathbf h,t} - s_t\|^2 +
 2\langle\widehat s_{\mathbf h,t} -\widehat s_{\mathbf h_0,t},
 s_t -\widehat s_{\mathbf h_0,t}\rangle -
 \|s_t -\widehat s_{\mathbf h_0,t}\|^2 +
 {\rm pen}(\mathbf h) - {\rm pen}(\widehat{\bf h}).
\end{eqnarray*}
Then,
\begin{eqnarray}
 \|\widehat s_{\widehat{\bf h},t} - s_t\|^2
 & \leqslant &
 \|\widehat s_{\mathbf h,t} - s_t\|^2 +
 2\langle\widehat s_{\mathbf h,t} -\widehat s_{\mathbf h_0,t},
 s_t -\widehat s_{\mathbf h_0,t}\rangle
 \nonumber\\
 & &
 \hspace{3cm} +
 {\rm pen}(\mathbf h) - {\rm pen}(\widehat{\bf h}) +
 2\langle\widehat s_{\widehat{\bf h},t} -\widehat s_{\mathbf h_0,t},
 \widehat s_{\mathbf h_0,t} - s_t\rangle
 \nonumber\\
 & = &
 \|\widehat s_{\mathbf h,t} - s_t\|^2 +
 {\rm pen}(\mathbf h) - {\rm pen}(\widehat{\bf h}) +
 2\langle\widehat s_{\widehat{\bf h},t} -\widehat s_{\mathbf h,t},
 \widehat s_{\mathbf h_0,t} - s_t\rangle
 \nonumber\\
 \label{risk_bound_numerator_PCO_NW_1}
 & = &
 \|\widehat s_{\mathbf h,t} - s_t\|^2 -\psi(\mathbf h) +\psi(\widehat{\bf h}),
\end{eqnarray}
where 
\begin{displaymath}
\psi(\cdot) :=
2\langle\widehat s_{.,t} - s_t,\widehat s_{\mathbf h_0,t} - s_t\rangle -
\textrm{pen}(\cdot).
\end{displaymath}
Now, let us rewrite $\psi(\cdot)$ in terms of $U_{.,\mathbf h_0}$, $W_{.,\mathbf h_0}$ and $W_{\mathbf h_0,.}$. For any $\mathbf h\in\mathcal H_{N}^{2}$,
\begin{eqnarray*}
 \psi(\mathbf h) & = &
 2\langle\widehat s_{\mathbf h,t} - s_{\mathbf h,t} + s_{\mathbf h,t} - s_t,
 \widehat s_{\mathbf h_0,t} - s_{\mathbf h_0,t} +
 s_{\mathbf h_0,t} - s_t\rangle - {\rm pen}(\mathbf h)\\
 & = &
 2\langle\widehat s_{\mathbf h,t} - s_{\mathbf h,t},
 \widehat s_{\mathbf h_0,t} - s_{\mathbf h_0,t}\rangle - {\rm pen}(\mathbf h)\\
 & &
 \hspace{2.5cm} +
 2\underbrace{(W_{\mathbf h,\mathbf h_0} +
 W_{\mathbf h_0,\mathbf h} +\langle s_{\mathbf h,t} - s_t,
 s_{\mathbf h_0,t} - s_t\rangle)}_{=:\psi_3(\mathbf h)}\\
 & = &
 \frac{2U_{\mathbf h,\mathbf h_0}}{N^2} +
 2\underbrace{\left(\frac{1}{N^2}\sum_{i = 1}^{N}
 \langle\Phi_{\mathbf h,t}(X^i;\cdot) - s_{\mathbf h,t},
 \Phi_{\mathbf h_0,t}(X^i;\cdot) - s_{\mathbf h_0,t}\rangle -
 \frac{{\rm pen}(\mathbf h)}{2}\right)}_{=:\psi_2(\mathbf h)} + 2\psi_3(\mathbf h)
\end{eqnarray*}
and, by the definition of ${\rm pen}(\mathbf h)$,
\begin{displaymath}
\psi_2(\mathbf h) = -\frac{1}{N^2}\left(\sum_{i = 1}^{N}
\langle\Phi_{\mathbf h,t}(X^i;\cdot),s_{\mathbf h_0,t}\rangle +
\sum_{i = 1}^{N}\langle\Phi_{\mathbf h_0,t}(X^i;\cdot),s_{\mathbf h,t}\rangle\right) +
\frac{1}{N}\langle s_{\mathbf h,t},s_{\mathbf h_0,t}\rangle.
\end{displaymath}
So,
\begin{displaymath}
\psi(\mathbf h) = 2(\psi_1(\mathbf h) + \psi_2(\mathbf h) + \psi_3(\mathbf h))
\quad {\rm with}\quad
\psi_1(\mathbf h) =\frac{U_{\mathbf h,\mathbf h_0}}{N^2}.
\end{displaymath}
\textbf{Step 2.} This step deals with suitable bounds on the $\psi_j$'s.
\begin{itemize}
 \item Consider $\mathbf h\in\mathcal H_N$. By Lemma \ref{control_U_term}, for any $\lambda > 0$ and $\theta\in (0,1)$, with probability larger than $1 - 5.4|\mathcal H_N|^2e^{-\lambda}$,
 \begin{eqnarray*}
  |\psi_1(\mathbf h)|
  & \leqslant &
  \frac{\theta}{2N}\overline s_{\mathbf h,t} +
  \frac{2\mathfrak c_{\ref{control_U_term}}(1 +\lambda)^3}{\theta N}\\
  & &
  \hspace{1.5cm}{\rm and}\quad
  |\psi_1(\widehat{\bf h})|
  \leqslant\frac{\theta}{2N}\overline s_{\widehat{\bf h},t} +
  \frac{2\mathfrak c_{\ref{control_U_term}}(1 +\lambda)^3}{\theta N}.
 \end{eqnarray*}
 \item For any $\mathbf h,\mathbf l\in\mathcal H_{N}^{2}$, consider 
 \begin{displaymath}
 \overline\psi_2(\mathbf h,\mathbf l) :=
 \frac{1}{N}\sum_{i = 1}^{N}
 \langle\Phi_{\mathbf h,t}(X^i;\cdot),s_{\mathbf l,t}\rangle.
 \end{displaymath}
 By Proposition \ref{Kernel_type_properties_Phi}.(4),
 \begin{displaymath}
 |\overline\psi_2(\mathbf h,\mathbf l)|\leqslant
 \mathfrak m_{\Phi}.
 \end{displaymath}
 Moreover, since $s_t\in\mathbb L^2(\mathbb R^2)$ by Inequality (\ref{square_integrability_s_t}),
 \begin{displaymath}
 |\langle s_{\mathbf h,t},s_{\mathbf h_0,t}\rangle|\leqslant
 \|Q_{\bf h}\star s_t\|\cdot\|Q_{\mathbf h_0}\star s_t\|\leqslant
 \|K\|_{1}^{4}\|s_t\|^2.
 \end{displaymath}
 Then, there exists a deterministic constant $\mathfrak c_1 > 0$, not depending on $N$ and $t$, such that
 \begin{displaymath}
 |\psi_2(\mathbf h)|\vee |\psi_2(\widehat{\bf h})|\leqslant
 \sup_{\mathbf l\in\mathcal H_{N}^{2}}
 |\psi_2(\mathbf l)|\leqslant
 \frac{\mathfrak c_1}{N}.
 \end{displaymath}
\item Consider $\mathbf h\in\mathcal H_{N}^{2}$. By Lemma \ref{control_W_term} and Cauchy-Schwarz's inequality, with probability larger than $1 - |\mathcal H_N|^2e^{-\lambda}$,
 \begin{eqnarray*}
  |\psi_3(\mathbf h)|
  & \leqslant &
  \frac{\theta}{4}(\|s_{\mathbf h,t} - s_t\|^2 +
  \|s_{\mathbf h_0,t} - s_t\|^2) +
  \frac{8\mathfrak c_{\ref{control_W_term}}(1 +\lambda)^2}{\theta N}\\
  & &
  \hspace{1.5cm}
  + 2\left(\frac{1}{2}\right)^{\frac{1}{2}}\left(\frac{\theta}{2}\right)^{\frac{1}{2}}
  \|s_{\mathbf h,t} - s_t\|\times
  \left(\frac{1}{2}\right)^{\frac{1}{2}}\left(\frac{2}{\theta}\right)^{\frac{1}{2}}
  \|s_{\mathbf h_0,t} - s_t\|\\
  & \leqslant &
  \frac{\theta}{2}\|s_{\mathbf h,t} - s_t\|^2 +
  \left(\frac{\theta}{4} +\frac{1}{\theta}\right)
  \|s_{\mathbf h_0,t} - s_t\|^2 +
  \frac{8\mathfrak c_{\ref{control_W_term}}(1 +\lambda)^2}{\theta N}
 \end{eqnarray*}	
 and
 \begin{displaymath}
 |\psi_3(\widehat{\bf h})|\leqslant
 \frac{\theta}{2}\|s_{\widehat{\bf h},t} - s_t\|^2 +
 \left(\frac{\theta}{4} +\frac{1}{\theta}\right)
 \|s_{\mathbf h_0,t} - s_t\|^2 +
 \frac{8\mathfrak c_{\ref{control_W_term}}(1 +\lambda)^2}{\theta N}.
 \end{displaymath}
\end{itemize}
\textbf{Step 3.} Let us establish that there exist two deterministic constants $\mathfrak c_2,\overline{\mathfrak c}_2 > 0$, not depending on $N$, $t$ and $\theta$, such that with probability larger than $1 -\overline{\mathfrak c}_2|\mathcal H_N|^2e^{-\lambda}$,
\begin{displaymath}
\sup_{\mathbf h\in\mathcal H_{N}^{2}}\left\{
\|\widehat s_{\mathbf h,t} - s_t\|^2 -
(1 +\theta)\left(\|s_{\mathbf h,t} - s_t\|^2 +
\frac{\overline s_{\mathbf h,t}}{N}\right)\right\}
\leqslant
\frac{\mathfrak c_2(1 +\lambda)^3}{\theta N}
\end{displaymath}
and
\begin{displaymath}
\sup_{\mathbf h\in\mathcal H_{N}^{2}}\left\{
\|s_{\mathbf h,t} - s_t\|^2 +\frac{\overline s_{\mathbf h,t}}{N} -
\frac{1}{1 -\theta}\|\widehat s_{\mathbf h,t} - s_t\|^2\right\}
\leqslant
\frac{\mathfrak c_2(1 +\lambda)^3}{\theta(1 -\theta)N}.
\end{displaymath}
On the one hand, for any $\mathbf h\in\mathcal H_{N}^{2}$,
\begin{eqnarray*}
 & &
 \|\widehat s_{\mathbf h,t} - s_t\|^2 -
 (1 +\theta)\left(\|s_{\mathbf h,t} - s_t\|^2 +\frac{\overline s_{\mathbf h,t}}{N}\right)\\
 & &
 \hspace{2cm} =
 \|\widehat s_{\mathbf h,t} - s_{\mathbf h,t}\|^2 +
 2\langle\widehat s_{\mathbf h,t} - s_{\mathbf h,t},
 s_{\mathbf h,t} - s_t\rangle +\|s_{\mathbf h,t} - s_t\|^2 -
 (1 +\theta)\left(\|s_{\mathbf h,t} - s_t\|^2 +\frac{\overline s_{\mathbf h,t}}{N}\right)\\
 & &
 \hspace{2cm} =
 \|\widehat s_{\mathbf h,t} - s_{\mathbf h,t}\|^2 -\frac{1 +\theta}{N}\overline s_{\mathbf h,t} +
 2W_{\mathbf h,\mathbf h} -\theta\|s_{\mathbf h,t} - s_t\|^2
\end{eqnarray*}
and
\begin{equation}\label{risk_bound_numerator_PCO_NW_2}
\|\widehat s_{\mathbf h,t} - s_{\mathbf h,t}\|^2 =
\frac{U_{\mathbf h,\mathbf h}}{N^2} +\frac{V_{\bf h}}{N}.
\end{equation}
So, with probability larger than $1 -\overline{\mathfrak c}_2|\mathcal H_N|^2e^{-\lambda}$,
\begin{eqnarray*}
 \sup_{\mathbf h\in\mathcal H_{N}^{2}}\left\{
 \|\widehat s_{\mathbf h,t} - s_{\mathbf h,t}\|^2 -
 \frac{1 +\theta}{N}\overline s_{\mathbf h,t}\right\}
 & \leqslant &
 \sup_{\mathbf h\in\mathcal H_{N}^{2}}\left\{
 \frac{|U_{\mathbf h,\mathbf h}|}{N^2} -\frac{\theta}{2N}\overline s_{\mathbf h,t} +
 \frac{1}{N}|V_{\bf h} -\overline s_{\mathbf h,t}| -
 \frac{\theta}{2N}\overline s_{\mathbf h,t}\right\}\\
 & \leqslant &
 \frac{2(\mathfrak c_{\ref{control_U_term}} +
 \mathfrak c_{\ref{control_V_term}})(1 +\lambda)^3}{\theta N}
\end{eqnarray*}
by Lemmas \ref{control_U_term} and \ref{control_V_term}, and then
\begin{equation}\label{risk_bound_numerator_PCO_NW_3}
\sup_{\mathbf h\in\mathcal H_{N}^{2}}\left\{
\|\widehat s_{\mathbf h,t} - s_t\|^2 -
(1 +\theta)\left(\|s_{\mathbf h,t} - s_t\|^2 +\frac{\overline s_{\mathbf h,t}}{N}\right)\right\}
\leqslant\frac{\mathfrak c_2(1 +\lambda)^3}{\theta N}
\end{equation}
by Lemma \ref{control_W_term}. On the other hand, for any $\mathbf h\in\mathcal H_{N}^{2}$,
\begin{eqnarray*}
 & &
 (1 -\theta)\left(\|s_{\mathbf h,t} - s_t\|^2 +\frac{\overline s_{\mathbf h,t}}{N}\right) -
 \|\widehat s_{\mathbf h,t} - s_t\|^2\\
 & &
 \hspace{2cm} =
 (1 -\theta)\left(\|s_{\mathbf h,t} - s_t\|^2 +\frac{\overline s_{\mathbf h,t}}{N}\right) -
 (\|\widehat s_{\mathbf h,t} - s_{\mathbf h,t}\|^2 + 2W_{\mathbf h,\mathbf h} +
 \|s_{\mathbf h,t} - s_t\|^2)\\
 & &
 \hspace{2cm} =
 -\theta\|s_{\mathbf h,t} - s_t\|^2 + (1 -\theta)\frac{\overline s_{\mathbf h,t}}{N}
 -\|\widehat s_{\mathbf h,t} - s_{\mathbf h,t}\|^2 - 2W_{\mathbf h,\mathbf h}\\
 & &
 \hspace{2cm}\leqslant
 2|W_{\mathbf h,\mathbf h}| -\theta\|s_{\mathbf h,t} - s_t\|^2 +
 \underbrace{\left|\frac{\overline s_{\mathbf h,t}}{N} -
 \|\widehat s_{\mathbf h,t} - s_{\mathbf h,t}\|^2\right|}_{=:\Lambda_{\bf h}} -
 \frac{\theta}{N}\overline s_{\mathbf h,t}
\end{eqnarray*}
and
\begin{displaymath}
 \Lambda_{\bf h} =
 \left|\frac{U_{\mathbf h,\mathbf h}}{N^2} +\frac{V_{\bf h}}{N} -
 \frac{\overline s_{\mathbf h,t}}{N}\right|
 \quad\textrm{by Equality (\ref{risk_bound_numerator_PCO_NW_2}).}
\end{displaymath}
By Lemmas \ref{control_U_term} and \ref{control_V_term}, there exist two deterministic constants $\mathfrak c_3,\overline{\mathfrak c}_3 > 0$, not depending $N$, $t$ and $\theta$, such that with probability larger than $1 -\overline{\mathfrak c}_3|\mathcal H_N|^2e^{-\lambda}$,
\begin{displaymath}
\sup_{\mathbf h\in\mathcal H_{N}^{2}}\left\{
\Lambda_{\bf h} -\frac{\theta}{N}\overline s_{\mathbf h,t}\right\}
\leqslant
\frac{\mathfrak c_3(1 +\lambda)^3}{\theta N}.
\end{displaymath}
Moreover, by Lemma \ref{control_W_term}, with probability larger than $1 - 2|\mathcal H_N|^2e^{-\lambda}$,
\begin{eqnarray*}
 \sup_{\mathbf h\in\mathcal H_{N}^{2}}\{2|W_{\mathbf h,\mathbf h}| -
 \theta\|s_{\mathbf h,t} - s_t\|^2\}
 & = &
 2\sup_{\mathbf h\in\mathcal H_{N}^{2}}\left\{
 |W_{\mathbf h,\mathbf h}| -\frac{\theta}{2}\|s_{\mathbf h,t} - s_t\|^2\right\}\\
 & \leqslant &
 \frac{4\mathfrak c_{\ref{control_W_term}}(1 +\lambda)^2}{\theta N}.
\end{eqnarray*}
So, with probability larger than $1 -\overline{\mathfrak c}_2|\mathcal H_N|^2e^{-\lambda}$,
\begin{equation}\label{risk_bound_numerator_PCO_NW_4}
\sup_{\mathbf h\in\mathcal H_{N}^{2}}\left\{
\|s_{\mathbf h,t} - s_t\|^2 +\frac{\overline s_{\mathbf h,t}}{N} -
\frac{1}{1 -\theta}\|\widehat s_{\mathbf h,t} - s_t\|^2\right\}
\leqslant
\frac{\mathfrak c_2(1 +\lambda)^3}{\theta(1 -\theta)N}.
\end{equation}
\textbf{Step 4.} By Step 2, there exist two deterministic constants $\mathfrak c_4,\overline{\mathfrak c}_4 > 0$, not depending on $N$, $t$ and $\theta$, such that with probability larger than $1 -\overline{\mathfrak c}_4|\mathcal H_N|^2e^{-\lambda}$,
\begin{displaymath}
|\psi(\mathbf h)|
\leqslant\theta\left(\|s_{\mathbf h,t} - s_t\|^2 +\frac{\overline s_{\mathbf h,t}}{N}\right) +
\left(\frac{\theta}{2} +\frac{2}{\theta}\right)\|s_{\mathbf h_0,t} - s_t\|^2 +
\frac{\mathfrak c_4(1 +\lambda)^3}{\theta N}
\end{displaymath} 
for every $\mathbf h\in\mathcal H_{N}^{2}$, and
\begin{displaymath}
|\psi(\widehat{\bf h})|
\leqslant\theta\left(\|s_{\widehat{\bf h},t} - s_t\|^2 +
\frac{\overline s_{\widehat{\bf h},t}}{N}\right) +
\left(\frac{\theta}{2} +\frac{2}{\theta}\right)\|s_{\mathbf h_0,t} - s_t\|^2 +
\frac{\mathfrak c_4(1 +\lambda)^3}{\theta N}.
\end{displaymath} 
So, by Inequality (\ref{risk_bound_numerator_PCO_NW_4}) (see Step 3), there exist two deterministic constants $\mathfrak c_5,\overline{\mathfrak c}_5 > 0$, not depending on $N$, $t$ and $\theta$, such that with probability larger than $1 -\overline{\mathfrak c}_5|\mathcal H_N|^2e^{-\lambda}$,
\begin{displaymath}
|\psi(\mathbf h)|
\leqslant
\frac{\theta}{1 -\theta}
\|\widehat s_{\mathbf h,t} - s_t\|^2 +
\left(\frac{\theta}{2} +\frac{2}{\theta}\right)\|s_{\mathbf h_0,t} - s_t\|^2 +
\mathfrak c_5\left(\frac{1}{\theta} +\frac{1}{1 -\theta}\right)
\frac{(1 +\lambda)^3}{N}
\end{displaymath}	
for every $\mathbf h\in\mathcal H_{N}^{2}$, and
\begin{displaymath}
|\psi(\widehat{\bf h})|
\leqslant
\frac{\theta}{1 -\theta}
\|\widehat s_{\widehat{\bf h},t} - s_t\|^2 +
\left(\frac{\theta}{2} +\frac{2}{\theta}\right)\|s_{\mathbf h_0,t} - s_t\|^2 +
\mathfrak c_5\left(\frac{1}{\theta} +\frac{1}{1 -\theta}\right)
\frac{(1 +\lambda)^3}{N}.
\end{displaymath}	
By Inequality \eqref{risk_bound_numerator_PCO_NW_1} (see Step 1), there exist two deterministic constants $\mathfrak c_6,\overline{\mathfrak c}_6 > 0$, not depending on $N$, $t$ and $\theta$, such that with probability larger than $1 -\overline{\mathfrak c}_6|\mathcal H_N|^2e^{-\lambda}$,
\begin{eqnarray*}
 \left(1 -\frac{\theta}{1 -\theta}\right)
 \|\widehat s_{\widehat{\bf h},t} - s_t\|^2
 & \leqslant &
 \left(1 +\frac{\theta}{1 -\theta}\right)
 \|\widehat s_{\mathbf h,t} - s_t\|^2\\
 & &
 \hspace{1.5cm} +
 \frac{\mathfrak c_6}{\theta}\left(\|s_{\mathbf h_0,t} - s_t\|^2 +
 \frac{(1 +\lambda)^3}{(1 -\theta)N}\right)
 \textrm{ $;$ }\forall\mathbf h\in\mathcal H_{N}^{2}.
\end{eqnarray*}	
By taking $\theta\in (0,1/2)$, the conclusion comes from Inequality (\ref{risk_bound_numerator_PCO_NW_3}) (see Step 3).\quad $\Box$
%


%
\subsubsection{Proof of Proposition \ref{Kernel_type_properties_Phi}}\label{proof_kernel_type_properties_Phi}
Let us establish the four {\it kernel type} properties of the map $\Phi_{\mathbf h,t}$, $\mathbf h\in (0,1]^2$, stated in Proposition \ref{Kernel_type_properties_Phi}. First of all, note that by Inequalities (\ref{properties_transition_density_1}) and (\ref{properties_transition_density_2}),
\begin{equation}\label{Kernel_type_properties_Phi_1}
\|s_t\|_{\infty} =
\sup_{(x,y)\in\mathbb R^2}|f(x)p_t(x,y)|\leqslant
\mathfrak c_1 :=\mathfrak m_p(t_0,T)\mathfrak m_f(t_0,T).
\end{equation}
\begin{enumerate}
 \item For every $\mathbf h = (h_1,h_2)\in (0,1]^2$ and $\varphi\in C^0([0,T])$, by Jensen's inequality and the change of variables formula,
 \begin{eqnarray*}
  \|\Phi_{\mathbf h,t}(\varphi;\cdot)\|^2
  & = &
  \int_{\mathbb R^2}\left(\frac{1}{T - t_0}\int_{t_0}^{T}Q_{\bf h}(\varphi(s) - x,\varphi(s + t) - y)ds\right)^2dxdy\\
  & \leqslant &
  \frac{1}{T - t_0}\int_{t_0}^{T}\underbrace{\left(
  \int_{-\infty}^{\infty}K_{h_1}(\varphi(s) - x)^2dx\right)}_{=\|K\|^2/h_1}
  \underbrace{\left(
  \int_{-\infty}^{\infty}K_{h_2}(\varphi(s + t) - y)^2dy\right)}_{=\|K\|^2/h_2}ds =
  \frac{\|K\|^4}{h_1h_2}.
 \end{eqnarray*}
 \item For any $\mathbf h = (h_1,h_2)$ and $\mathbf l = (\ell_1,\ell_2)$ belonging to $(0,1]^2$, by Jensen's inequality (three times), and since $X^1$ and $X^2$ are independent processes,
 \begin{eqnarray*}
  & &
  \mathbb E(\langle\Phi_{\mathbf h,t}(X^1;\cdot),
  \Phi_{\mathbf l,t}(X^2;\cdot)\rangle^2)\\
  & &
  \hspace{1.5cm} =
  \mathbb E\left[\left(\frac{1}{T - t_0}\int_{t_0}^{T}
  \int_{-\infty}^{\infty}K_{h_1}(X_{s}^{1} - x)
  \int_{-\infty}^{\infty}K_{h_2}(X_{s + t}^{1} - y)
  \Phi_{\mathbf l,t}(X^2;x,y)dydxds\right)^2
  \right]\\
  & &
  \hspace{1.5cm}\leqslant
  \|K\|_{1}^{2}
  \int_{\mathbb R^2}\left(\frac{1}{T - t_0}
  \int_{t_0}^{T}\mathbb E(|K_{h_1}(X_s - x)K_{h_2}(X_{s + t} - y)|)ds\right)
  \mathbb E(\Phi_{\mathbf l,t}(X;x,y)^2)dydx.
 \end{eqnarray*}
 Moreover, for every $(x,y)\in\mathbb R^2$, by Equality (\ref{consequence_homogeneity}) and Inequality (\ref{Kernel_type_properties_Phi_1}),
 \begin{eqnarray*}
  & &
  \frac{1}{T - t_0}\int_{t_0}^{T}
  \mathbb E(|K_{h_1}(X_s - x)K_{h_2}(X_{s + t} - y)|)ds\\
  & &
  \hspace{2cm} =
  \int_{\mathbb R^2}|K_{h_1}(\xi - x)K_{h_2}(\zeta - y)|\left(
  \frac{1}{T - t_0}\int_{t_0}^{T}p_{s,s + t}(\xi,\zeta)ds\right)d\xi d\zeta\\
  & &
  \hspace{2cm} =
  \int_{\mathbb R^2}|K_{h_1}(\xi - x)K_{h_2}(\zeta - y)|
  s_t(\xi,\zeta)d\xi d\zeta
  \leqslant
  \mathfrak c_1\|K\|_{1}^{2}.
 \end{eqnarray*}
 Therefore,
 \begin{displaymath}
 \mathbb E(\langle\Phi_{\mathbf h,t}(X^1;\cdot),
 \Phi_{\mathbf l,t}(X^2;\cdot)\rangle^2)
 \leqslant
 \mathfrak c_1\|K\|_{1}^{4}
 \underbrace{\int_{\mathbb R^2}\mathbb E(
 \Phi_{\mathbf l,t}(X;x,y)^2)dydx}_{=\overline s_{\mathbf l,t}}.
 \end{displaymath}
 \item For every $\mathbf h\in (0,1]^2$ and $\varphi\in\mathbb L^2(\mathbb R^2)$, by Jensen's inequality, Equality (\ref{consequence_homogeneity}) and Inequality (\ref{Kernel_type_properties_Phi_1}),
 \begin{eqnarray*}
  \mathbb E(\langle\Phi_{\mathbf h,t}(X;\cdot),\varphi\rangle^2) & = &
  \mathbb E\left[\left(\frac{1}{T - t_0}\int_{t_0}^{T}\int_{\mathbb R^2}
  Q_{\bf h}(X_s - x,X_{s + t} - y)\varphi(x,y)dxdyds\right)^2
  \right]\\
  & \leqslant &
  \int_{\mathbb R^2}(Q_{\bf h}\star\varphi)(\xi,\zeta)^2s_t(\xi,\zeta)d\xi d\zeta
  \leqslant\mathfrak c_1\|Q_{\bf h}\star\varphi\|^2
  \leqslant\mathfrak c_1\|K\|_{1}^{4}\|\varphi\|^2.
 \end{eqnarray*}
 \item For every $\mathbf h,\mathbf l\in (0,1]^2$, by Inequality (\ref{Kernel_type_properties_Phi_1}),
 \begin{eqnarray*}
  |\langle\Phi_{\mathbf h,t}(X,\cdot),s_{\mathbf l,t}\rangle| & = &
  \frac{1}{T - t_0}\left|\int_{t_0}^{T}(Q_{\bf h}\star s_{\mathbf l,t})(X_s,X_{s + t})ds\right|\\
  & \leqslant &
  \|Q_{\bf h}\star s_{\mathbf l,t}\|_{\infty}\leqslant
  \|Q_{\bf h}\|_1\|Q_{\bf l}\|_1\|s_t\|_{\infty}\leqslant
  \mathfrak c_1\|K\|_{1}^{4}.
  \quad\Box
 \end{eqnarray*}
\end{enumerate}
%


%
\subsection{Proof of Corollary \ref{risk_bound_PCO_NW}}
The proof of Corollary \ref{risk_bound_PCO_NW} relies on Theorem \ref{risk_bound_numerator_PCO_NW} and on the following technical lemma.
%


%
\begin{lemma}\label{concentration_inequalities_expectation}
Let $R$ be a random variable, and assume that there exist $r,c > 0$ and $q\geqslant 1$ such that, for every $\alpha\in\mathbb R_+$,
\begin{displaymath}
\mathbb P\left(R\leqslant\frac{\alpha^q}{r}\right)
\geqslant 1 - ce^{-\alpha}.
\end{displaymath}
Then,
\begin{displaymath}
\mathbb E(R)\leqslant
\frac{2^{q + 1}\log(c)^q}{r} +
\frac{\mathfrak c_q}{r}
\quad\textrm{with}\quad
\mathfrak c_q =
\int_{0}^{\infty}\exp\left(-\frac{1}{2}\beta^{\frac{1}{q}}\right)d\beta <\infty.
\end{displaymath}
\end{lemma}
The proof of Lemma \ref{concentration_inequalities_expectation} is postponed to Section \ref{proof_concentration_inequalities_expectation}.
%


%
\subsubsection{Steps of the proof of Corollary \ref{risk_bound_PCO_NW}}
As in the proof of Corollary \ref{risk_bound_NW}, there exists a constant $\mathfrak c_1 > 0$, not depending on $N$, $t$, $\texttt l$ and $\texttt r$, such that
\begin{eqnarray*}
 \mathbb E(\|\widehat p_{\widehat{\bf h},\widehat\ell,t} -
 p_t\|_{[\texttt l,\texttt r]\times\mathbb R}^{2})
 & \leqslant &
 \frac{\mathfrak c_1}{m^2}
 (\mathbb E(\|\widehat s_{\widehat{\bf h},t} - s_t\|^2) +
 \mathbb E(\|\widehat f_{\widehat\ell} - f\|^2)).
\end{eqnarray*}
Moreover, by Theorem \ref{risk_bound_numerator_PCO_NW}, Marie and Rosier (2023), Theorem 1, and by Lemma \ref{concentration_inequalities_expectation}, there exists a constant $\mathfrak c_2 > 0$, not depending on $N$, $t$, $\texttt l$ and $\texttt r$, such that
\begin{displaymath}
\mathbb E(\|\widehat s_{\widehat{\bf h},t} - s_t\|^2)
\leqslant
\mathfrak c_2\left(\min_{\mathbf h\in\mathcal H_{N}^{2}}
\mathbb E(\|\widehat s_{\mathbf h,t} - s_t\|^2) +
\|s_{\mathbf h_0,t} - s_t\|^2 +\frac{\log(N)^6}{N}\right)
\end{displaymath}
and
\begin{displaymath}
\mathbb E(\|\widehat f_{\widehat\ell} - f\|^2)
\leqslant
\mathfrak c_2\left(\min_{\ell\in\mathcal H_N}
\mathbb E(\|\widehat f_{\ell} - f\|^2) +
\|f_{h_0} - f\|^2 +\frac{\log(N)^3}{N}\right).
\end{displaymath}
Therefore, there exists a constant $\mathfrak c_3 > 0$, not depending on $N$, $t$, $\texttt l$ and $\texttt r$, such that
\begin{eqnarray*}
 \mathbb E(\|\widehat p_{\widehat{\bf h},\widehat\ell,t} -
 p_t\|_{[\texttt l,\texttt r]\times\mathbb R}^{2})
 & \leqslant &
 \frac{\mathfrak c_3}{m^2}\left(
 \min_{(\mathbf h,\ell)\in\mathcal H_{N}^{3}}
 \{\mathbb E(\|\widehat s_{\mathbf h,t} - s_t\|^2) +
 \mathbb E(\|\widehat f_{\ell} - f\|^2)\}\right.\\
 & &
 \hspace{3cm}\left. +
 \|s_{\mathbf h_0,t} - s_t\|^2 +\|f_{h_0} - f\|^2 +\frac{\log(N)^6}{N}\right)
 \quad\Box
\end{eqnarray*}
%


%
\subsubsection{Proof of Lemma \ref{concentration_inequalities_expectation}}\label{proof_concentration_inequalities_expectation}
Consider $A > 0$. First, by the Fubini-Tonelli theorem,
\begin{eqnarray*}
 \mathbb E(R) & = &
 \mathbb E(R\mathbf 1_{R\leqslant A}) +
 \mathbb E((R - A)\mathbf 1_{R > A}) +
 A\mathbb P(R > A)\\
 & \leqslant &
 2A +\mathbb E\left(\mathbf 1_{R > A}
 \int_{A}^{\infty}\mathbf 1_{R > x}dx\right)
 \leqslant
 2A +\int_{A}^{\infty}\mathbb P(R > x)dx.
\end{eqnarray*}
Now, by the change of variables formula,
\begin{eqnarray*}
 \int_{A}^{\infty}\mathbb P(R > x)dx
 & = &
 \frac{1}{r}\int_{rA}^{\infty}\mathbb P\left(R >\frac{\beta}{r}\right)d\beta\\
 & \leqslant &
 \frac{c}{r}\int_{rA}^{\infty}e^{-\beta^{1/q}}d\beta\
 \leqslant
 \frac{\mathfrak c_qc}{r}\exp\left(-\frac{1}{2}(rA)^{\frac{1}{q}}\right).
\end{eqnarray*}
Then, for $A =\log(c^2)^q/r$,
\begin{eqnarray*}
 \mathbb E(R) & \leqslant &
 \frac{2\log(c^2)^q}{r} +
 \frac{\mathfrak c_qc}{r}\exp\left(-\frac{1}{2}\log(c^2)\right)\\
 & = &
 \frac{2^{q + 1}\log(c)^q}{r} +
 \frac{\mathfrak c_q}{r}
 \quad\Box
\end{eqnarray*}
%


%

%

\begin{thebibliography}{99}
 \bibitem{AS99}	 A\"it-Sahalia, Y. (1999). Transition Densities for Interest Rate and Other Nonlinear Diffusions. {\it The Journal of Finance} LIV, 1361-1395.
 \bibitem{COMTE17} Comte, F. (2017). {\it Estimation non-param\'etrique.} Paris: Spartacus IDH.
 \bibitem{CGC20} Comte, F. and Genon-Catalot, V. (2020). Nonparametric Drift Estimation for I.I.D. Paths of Stochastic Differential Equations. {\it The Annals of Statistics} 48, 6, 3336-3365.
 \bibitem{CGC23} Comte, F. and Genon-Catalot, V. (2024). Nonparametric Estimation for I.I.D. Stochastic Differential Equations with Space-Time Dependent Coefficients. {\it SIAM/ASA Journal of Uncertainty Quantification} 12, 2, 377-410.
 \bibitem{CM21} Comte, F. and Marie, N. (2021). Nonparametric Estimation for I.I.D. Paths of Fractional SDE. {\it Statistical Inference for Stochastic Processes} 24, 3, 669-705.
 \bibitem{CM22} Comte, F. and Marie, N. (2023). Nonparametric Drift Estimation from Diffusions with Correlated Brownian Motions. {\it Journal of Multivariate Analysis} 198, 23 pages.
 \bibitem{CM25} Comte, F. and Marie, N. (2025). Nonparametric Estimation of the Transition Density Function for Diffusion Processes. {\it Stochastic Processes and their Applications} 188, 27 pages.
 \bibitem{DMH22} Della Maestra, L. and  Hoffmann, M. (2022). Nonparametric Estimation for Interacting Particle Systems: McKean-Vlasov Models. {\it Probability Theory Related Fields} 182, 551-613.
 \bibitem{DDBM21} Denis, C., Dion-Blanc, C. and Martinez, M. (2021). A Ridge Estimator of the Drift from Discrete Repeated Observations of the Solution of a Stochastic Differential Equation. {\it Bernoulli} 27, 2675-2713.
 \bibitem{GN15} Gin\'e, E. and Nickl, R. (2015). {\it Mathematical Foundations of Infinite-Dimensional Statistical Models.} Cambridge: Cambridge University Press.
 \bibitem{HM20} Halconruy, H. and Marie, N. (2020). Kernel Selection in Nonparametric Regression. {\it Mathematical Methods of Statistics} 29, 32-56.
 \bibitem{HM23} Halconruy, H. and Marie, N. (2023). On a Projection Least Squares Estimator for Jump Diffusion Processes. {\it Annals of the Institute of Statistical Mathematics} 76, 2, 209-234.
 \bibitem{KS85} Kusuoka, S. and Stroock, D. (1985). Applications of the Malliavin Calculus, Part II. {\it Journal of the Faculty of Science, University of Tokyo} 32, 1-76.
 \bibitem{KUTOYANTS04} Kutoyants, Y. (2004). {\it Statistical Inference for Ergodic Diffusion Processes.} London: Springer.
 \bibitem{LACOUR07} Lacour, C. (2007). Adaptive Estimation of the Transition Density of a Markov Chain. {\it Annales de l'Institut Henri Poincar\'e B} 43, 571-597.
 \bibitem{LACOUR08} Lacour, C. (2008) Nonparametric Estimation of the Stationary Density and the Transition Density of a Markov Chain. {\it Stochastic Processes and their Applications} 118, 232-260.
 \bibitem{LMR17} Lacour, C., Massart, P. and Rivoirard, V. (2017). Estimator Selection: a New Method with Applications to Kernel Density Estimation. {\it Sankhya} 79, 298-335.
 \bibitem{MARIE23} Marie, N. (2023). Nonparametric Estimation for I.I.D. Paths of a Martingale Driven Model with Application to Non-Autonomous Financial Models. {\it Finance and Stochastics} 27, 1, 97-126.
 \bibitem{MR23} Marie, N. and Rosier, A. (2023). Nadaraya-Watson Estimator for I.I.D. Paths of Diffusion Processes. {\it Scandinavian Journal of Statistics} 50, 2, 589-637.
 \bibitem{MASSART07} Massart, P. (2007). {\it Concentration Inequalities and Model Selection.} Berlin, Heidelberg: Springer.
 \bibitem{MPZ21} Menozzi, S., Pesce, A. and Zhang, X. (2021). Density and Gradient Estimates for Non Degenerate Brownian SDEs with Unbounded Measurable Drift. {\it Journal of Differential Equations} 272, 330-369.
 \bibitem{MSS04} Milstein, G.N., Schoenmakers, J.G.M. and Spokoiny, V. (2004). Transition Density Estimation for Stochastic Differential Equations via Forward-Reverse Representations. {\it Bernoulli} 10(2), 281-312.
 \bibitem{RS07} Ramsay, J.O. and Silverman, B.W. (2007). {\it Applied Functional Data Analysis: Methods and Case Studies.} New-York: Springer.
 \bibitem{SART14} Sart, M. (2014). Estimation of the Transition Density of a Markov Chain. {\it Annales de l'Institut Henri Poincar\'e B} 50, 1028-1068.
 \bibitem{WCM16} Wang, J.-L., Chiou, J.-M. and Mueller, H.-G. (2016). Review of Functional Data Analysis. {\it Annual Review of Statistics and its Application} 3, 257-295.
\end{thebibliography}
\end{document}